\documentclass[journalpaper, 9pt]{article}
\usepackage[paper=letterpaper,centering,margin=1.0in]{geometry}

\usepackage[american]{babel}
\usepackage{epsfig}
\usepackage{color}
\usepackage{amsthm,amsfonts,amsmath,amssymb,amstext,latexsym}

\usepackage{color}
 \usepackage{algorithm}
\usepackage{algpseudocode}

\allowdisplaybreaks

\newtheorem{theorem}{Theorem}
\newtheorem{lemma}{Lemma}

\newtheorem{remark}{Remark}

\newtheorem{proposition}{Proposition}

\newtheorem{definition}{Definition}
\newtheorem{corollary}{Corollary}

\newcommand{\bydef}{\stackrel{\rm{def}}{=}}

\title{Power-of-d-Choices with Memory: Fluid Limit and Optimality}

\author{
Jonatha Anselmi\footnote{INRIA Bordeaux Sud Ouest, 200 av. de la Vieille Tour, 33405 Talence,  France. Email: jonatha.anselmi@inria.fr} ~and Francois Dufour\footnote{INRIA Bordeaux Sud Ouest, 200 av. de la Vieille Tour, 33405 Talence,  France. Email: francois.dufour@math.u-bordeaux.fr}
}

\begin{document}
\date{}
\maketitle

\begin{abstract}
In multi-server distributed queueing systems, the access of stochastically arriving jobs to resources is often regulated by a dispatcher, also known as load balancer. A fundamental problem consists in designing a load balancing algorithm that minimizes the delays experienced by jobs. During the last two decades, the power-of-$d$-choice algorithm, based on the idea of dispatching each job to the least loaded server out of $d$ servers randomly sampled at the arrival of the job itself, has emerged as a breakthrough in the foundations of this area due to its versatility and appealing asymptotic properties. In this paper, we consider the power-of-$d$-choice algorithm with the addition of a local memory that keeps track of the latest observations collected over time on the sampled servers. Then, each job is sent to a server with the lowest observation. We show that this algorithm is asymptotically optimal in the sense that the load balancer can always assign each job to an idle server in the large-system limit. This holds true if and only if the system load~$\lambda$ is less than~$1-\frac{1}{d}$. If this condition is not satisfied, we show that queue lengths are tightly bounded by $\left\lceil - \frac{ \log (1-\lambda)}{\log (\lambda d +1)} \right\rceil$. This is in contrast with the classic version of the power-of-$d$-choice algorithm, where at the fluid scale a strictly positive proportion of servers containing $i$ jobs exists for all $i\ge 0$, in equilibrium. Our results quantify and highlight the importance of using memory as a means to enhance performance in randomized load balancing.
\end{abstract}


\section{Introduction}

\begin{figure}
\centering
\includegraphics[width=11cm]{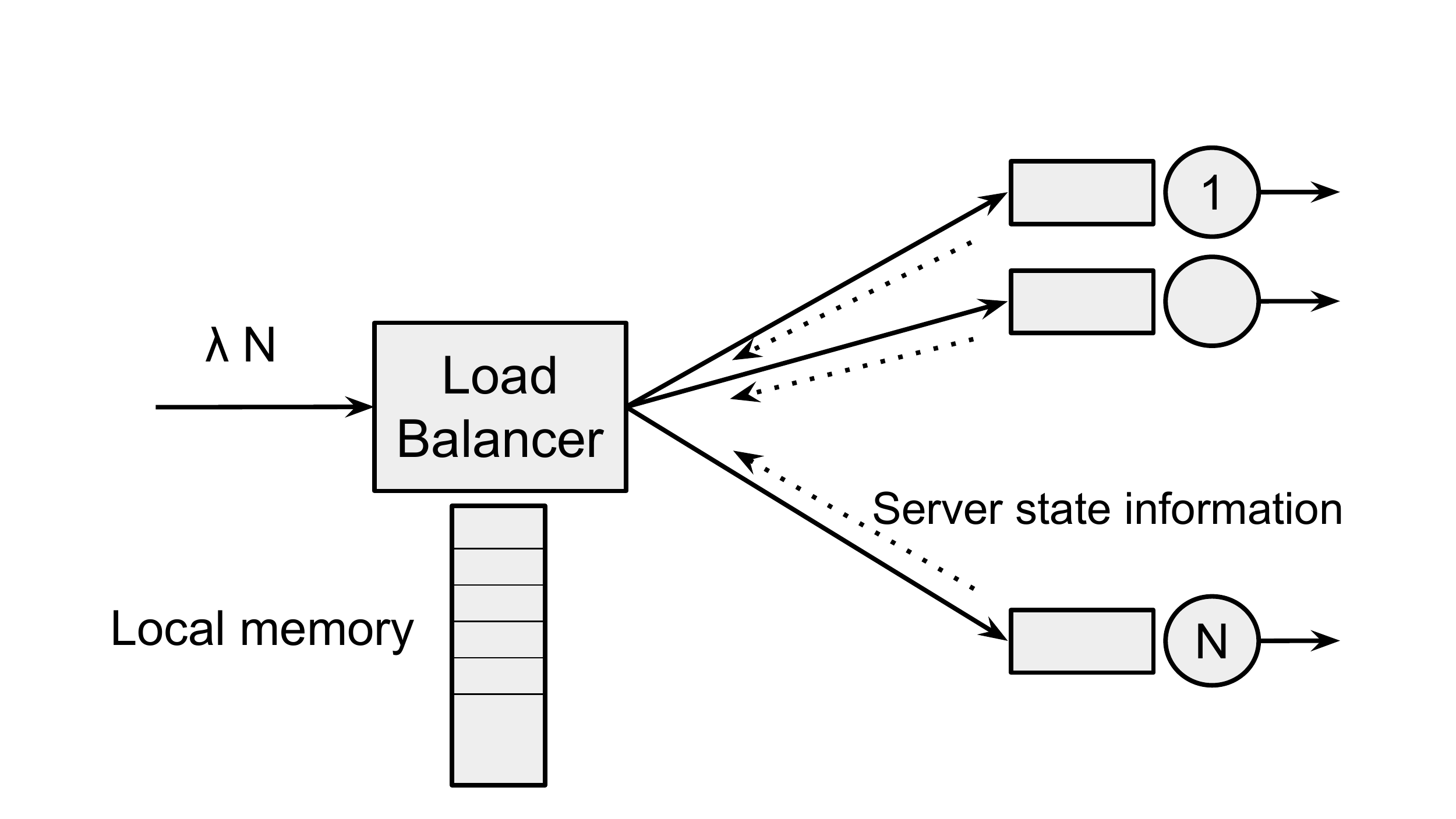}
\caption{Architecture of the distributed system for load balancing.}
\label{system}
\end{figure}

In multi-server distributed queueing systems, the access of stochastically arriving jobs to resources, or servers, is often regulated by a central dispatcher, also known as load balancer.
A fundamental problem consists in designing a load balancing algorithm able to minimize the delays experienced by jobs.
In this paper, we are interested in a setting where a traffic of rate $\lambda N$ needs to be distributed across~$N$ unit-rate parallel servers, each with its own queue, as indicated in Figure~\ref{system}.
The load balancer may rely on feedback information coming from the servers, which may also be stored in a local memory.
Depending on the architecture, feedback information can arrive at the dispatcher through a push- or pull-based mechanism.
In the former, the dispatcher initiates the communication fetching the requested information from the servers, while in the latter servers periodically send state information to the dispatcher.
%
This type of model finds applications in computer and communication systems, hospitals and road networks,
and there exists a significant and growing number of references;
see, e.g., the recent works
Ying et al. \cite{P2C_MOR17},
Gardner et al. \cite{RedundancyOR17},
Gamarnik et al. \cite{tsitsiklis2017},
Gupta and Walton \cite{2017Walton}
and the references therein.
Nevertheless, it is often difficult to establish whether an algorithm is better than another because in general the answer strongly depends on the underlying architecture, application or traffic conditions.
For instance, assigning jobs to servers uniformly at random or in a cyclic fashion provides a very scalable dispatching scheme as it requires neither static nor dynamic information about servers but the resulting performance is quite poor; the join-the-shortest-queue algorithm is ``optimal'' under some conditions, Winston \cite{winston1977} and Weber \cite{weber1978}, but its applicability in large systems is debated due to the high communication overhead between the servers and the dispatcher; the join-the-idle-queue algorithm, Lu et al. \cite{Lu2011}, performs very well when $N$ is large, see Stolyar~\cite{Stolyar2015}, but poorly when $\lambda$ gets close to its critical value, and in addition it requires servers to generate messages on their own; for a more complete discussion, we point the reader to the recent survey in Van der Boor et al. \cite{2018arXiv180605444V}.

During the last two decades, the power-of-$d$-choice algorithm, introduced in Mitzenmacher \cite{Mitzenmacher2001}, Vvedenskaya et al. \cite{Vvedenskaya} and referred to as SQ($d$), has emerged as a breakthrough in the foundations of this area due to its versatility and its appealing asymptotic properties.
It works as follows: upon arrival of each job, $d\ge 2$ servers are contacted uniformly at random, their state (e.g., queue length or workload) is retrieved, and then the job is dispatched to a server in the best observed state (among the $d$ selected).
The first remarkable property is that in the large-system limit, $N\to\infty$, the stationary proportion of servers with at least~$i$ jobs decreases doubly exponentially in~$i$, though it remains strictly positive for all $i$.
This result has been generalized in Bramson et al. \cite{Bramson2010} to the case where service times are heavy-tailed rather than exponential; see also Bramson et al. \cite{Bramson2012}.
In addition, it turns out that SQ($d$) is \emph{heavy-traffic optimal} in the sense that it minimizes the workload or queue-length process over all time in the diffusion limit where $\lambda \uparrow 1$; see Chen and Ye \cite{Chen2012} and Maguluri et al. \cite{Maguluri2012}.
In Ying et al. \cite{P2C_MOR17}, it is also shown that the number of sampled servers can be dramatically reduced if tasks arrive in batches, which is useful to reduce the communication overhead between the load balancer and the servers.
In Mitzenmacher et al. \cite{1182005}, the power-of-$d$-choice algorithm is studied in the case where the load balancer is endowed with a local memory that stores the index and the state of the least loaded server out of the $d$ sampled each time a job arrives.
When the $n$-th job arrives, the winning server is chosen among the $d$ servers randomly selected upon its arrival and the server associated to the observation stored in the memory.
The resulting performance is better than the one achieved by SQ$(2d)$.
In the standard memoryless case, if $d$ is allowed to depend on $N$ and $d(N)\to\infty$, SQ($d$) has been recently shown 
to become \emph{fluid (or mean-field) optimal}, i.e., optimal in the large-system limit, with a diffusion limit matching the one of the celebrated join-the-shortest-queue algorithm provided that $d(N)$ grows to infinity sufficiently fast; see Mukherjee et al. \cite{Borst2016}, Dieker and Suk \cite{dieker2015}.
At a fluid scale, optimality here is related to the ability of assigning each incoming job to an idle server.
Also our work aims at achieving fluid optimality but we will consider~$d$ as a constant to keep the communication overhead at a minimum.
Towards this purpose, we will show that it is enough to endow the load balancer with a local memory that keeps track of the latest observation collected on each server.
This approach is also close to Mitzenmacher et al. \cite{1182005}, though different because in that reference the memory can only store one observation.
In fact, one observation (or even a finite number of observations) is not enough to achieve fluid optimality; see Gamarnik et al.~\cite{Gamarnik2016}.
We observe that fluid optimality can also be achieved by the join-the-shortest-queue and join-the-idle-queue algorithms.
However, these are not directly comparable to our algorithm because they are meant to run on a different architecture (pull-based rather than push-pased).

The fact that we consider a memory with~$N$ slots has an impact on our proofs.
As discussed in Mitzenmacher et al. \cite{1182005}, if the memory size is uniformly bounded then the observations in the local memory evolve much faster than the actual queue lengths, and in this case to establish fluid limit results one can adopt the ad-hoc proof technique developed in Luczak and Norris~\cite{luczak2013}.
On the other hand, this does not apply to our case because observations and queue lengths evolve within the same timescale.
Also the pull-based version of join-the-idle-queue, Lu et al. \cite{Lu2011}, requires a memory with $N$ slots but the main difference with respect to our approach is that the information stored in the memory is always up to date, which is not the case within our algorithm.

\subsection{Contribution.}
In Algorithm~\ref{alg}, we provide a pseudocode for the proposed power-of-$d$-choices algorithm with memory and~$N$ servers, referred to as SQ$(d,N)$; some variants of such algorithm are also discussed in the Conclusions.
Upon arrival of one job, the states collected from~$d$ randomly chosen servers are stored in the local array \texttt{Memory}. Then, the job is sent to a server chosen randomly (with replacement) among the ones having the lowest recorded state.
Finally, the observation of the selected server is incremented by one.
\begin{algorithm}
\caption{Power-of-$d$-choices with memory and $N$ servers.}\label{alg}
\begin{algorithmic}[1]
\Procedure{SQ}{$d,N$}
\State Memory$[i]= 0,\,\,\forall i=1,\ldots,N$;
\For{each job arrival }
  \For{$i= 1,\ldots, d$}
    \State     rnd\_server = random(1,$\ldots,N$);
    \State     Memory[rnd\_server] = get\_state(rnd\_server);
  \EndFor
  \State  selected\_server = random($ \arg\min_{i\in\{1,\ldots,N\}}$ Memory[$i$]);
  \State  send\_job\_to(selected\_server);\
  \State  Memory[selected\_server]++;
\EndFor
\EndProcedure
\end{algorithmic}
\end{algorithm}

It is intuitive that SQ($d,N$) results in more balanced allocations than SQ$(d)$.
This follows by using the coupling argument developed in~Theorem~3.5 of Azar et al. \cite{Azar2000}, which can be adapted to argue that at any point in time the vector of queue lengths achieved with SQ$(d,N)$ is majorized by the vector of queue lengths achieved with SQ$(d)$.
On the other hand, it is not clear how much such improvement can be. This is the goal of the present paper.


We investigate the time-varying dynamics of SQ$(d,N)$ by means of a continuous-time Markov chain $X^N(t)$ that keeps track of the proportion of servers with $i$ jobs and for which their last observation collected by the load balancer is~$j$, for all $i$ and $j$.
To the best of our knowledge, this is the first paper that studies the dynamics induced by SQ$(d,N)$. 
The transition rates of $X^N(t)$ are non-Lipschitz and a satisfactory analysis of~$X^N(t)$ when~$N$ is finite seems to be out of reach.
Our main contributions are as follows:
\begin{enumerate}
\item In Theorem~\ref{th0}, we let $N\to\infty$ and identify the fluid limit of $X^N(t)$, an absolutely continuous function that is interpreted as a first-order approximation of the original model $X^N(t)$.
The fluid limit is motivated by the fact that real systems are composed of many servers and that it enables a tractable analysis for the dynamics of SQ$(d,N)$.
A fluid limit is necessarily a \emph{fluid solution}, as introduced in Definition~\ref{def1}.
The proof of the fluid limit is the main technical part of this work and is given in Section~\ref{sec:conn}.
The main difficulty stands in the discontinuous structure of the drift of $X^N(t)$; see Section~\ref{eq:cus78dc} for further details.
We obtain the fluid limit under a finite buffer assumption, though as discussed in the Conclusions we believe that this assumption can be relaxed.

\item We then study fixed points, fluid solutions that are constant over time.
Theorem~\ref{th1} shows that there exists a unique fixed point.
The general structure of such fixed point as a function of $\lambda$ is quite involved and implies that in equilibrium
\begin{itemize}
 \item [a)]
Fluid queue lengths are uniformly and tightly bounded by $j^\star+1$, where
\begin{equation}
\label{eq:jstar}
j^\star\bydef \left\lfloor   -\frac{ \log (1-\lambda)}{\log (\lambda d +1)} \right\rfloor.
\end{equation}
This is in contrast with SQ$(d)$, where queue lengths are unbounded in the sense that a strictly positive proportion of servers containing $i$ jobs exists for all $i\ge 0$, in the fluid equilibrium; see Mitzenmacher \cite{Mitzenmacher2001}.
Figure~\ref{fig:j_star} illustrates the behavior of the upper bound $j^\star+1$ by varying $\lambda$ and~$d$, and shows that the size of the most loaded server will remain very small even when $\lambda$ is very close to its critical value.
\begin{figure}
\centering
\includegraphics[width=11cm]{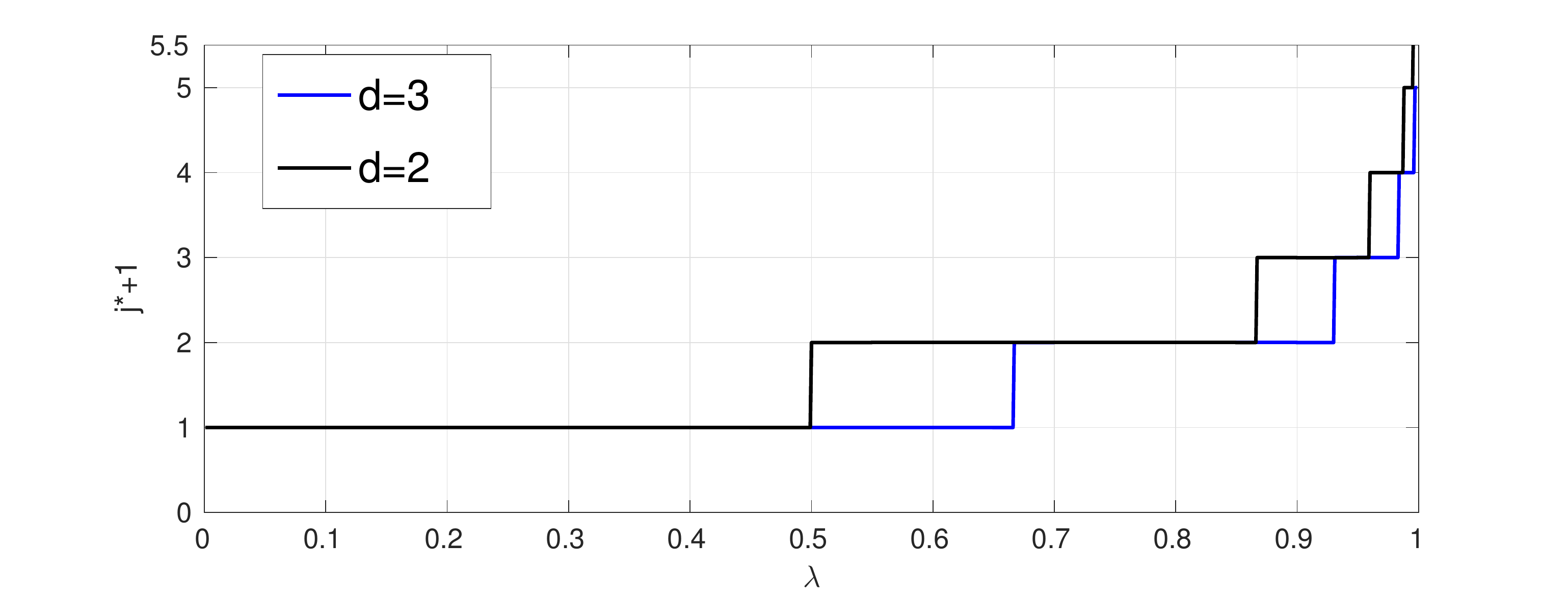}
\caption{Plots of the maximum queue length, $j^\star+1$, by varying $\lambda$ and~$d$.}
\label{fig:j_star}
\end{figure}
In fact, even when $\lambda = 0.995$ and $d=2$, at the fluid scale no server will contain more than just 5 jobs.

\item [b)]
The load balancer memory can only contain two possible observations, namely~$j^\star$ and $j^\star+1$.
\end{itemize}
The case of particular interest is when $\lambda < 1-1/d$, where $j^\star=0$ and thus the load balancer memory always contains a strictly positive proportion of zeros. This means that the load balancer can always assign incoming jobs to idle servers, which is clearly the ideal situation for any incoming job. In this sense we say that SQ$(d,N)$ is \emph{asymptotically optimal}.
When $\lambda \ge 1-1/d$ the load balancer memory will never contain a strictly positive mass of zeros but it will still be able to assign a fraction of jobs to idle servers ensuring that the average number of jobs in each queue belongs to the interval   $j^\star -\frac{1}{d} + \frac{1}{2} \pm \frac{1}{2}$ (Proposition~\ref{prop:LS_LM}).

\item Finally, we investigate stability properties of the unique fixed point.
Theorem~\ref{th2} establishes that fluid solutions converge to such point regardless of the initial condition and exponentially fast, provided that $\lambda<1-1/d$.
Thus, in this case all fluid solutions will be eventually asymptotically optimal as the load balancer memory will eventually be populated by a strictly positive mass of zeros.
The proof of this result, given in Section~\ref{proof_th2}, is based on a sort of Lyapunov argument that allows us to show that the time evolution of fluid solutions is eventually governed by the unique solution of a linear and autonomous ODE system.

\end{enumerate}

In summary, the proposed algorithm SQ($d,N$) has the same communication overhead of its memoryless counterpart SQ($d$) but a much better performance, which is paid at the cost of endowing the controller with a memory of~$N$ slots.
It is to be noted that asymptotic optimality can be obtained

%
%
%
%

\section{Performance models}
\label{model}

In order to describe the time varying effects of SQ$(d,N)$ on queue lengths, we introduce a stochastic and a deterministic model.
The stochastic model is meant to capture the variability of job interarrival and service times that is intrinsic in 
multi-server distributed queueing systems.
Due to its intractability, a satisfactory analysis of such model is out of reach.
In this respect, the deterministic model is convenient because it does enable analytical tractability.
In this section, we also show our first result, which states that both models are connected each other: the deterministic can be interpreted as a first-order approximation of the stochastic.

In the following, we will refer to a server with~$i$ jobs and for which its last observation at the controller is~$j$ as an \emph{$(i,j)$-server}.

\subsection{Markov model.}

First, we model the dynamics induced by SQ($d,N$) as a Markov chain in continuous-time: arrivals at the load balancer are assumed to follow a Poisson process with rate~$\lambda N$, with $0<\lambda<1$, and service times are independent, exponentially distributed random variables with unit mean.
Servers process jobs according to any work-conserving discipline and each of them can contain $I>1$ jobs at most.
A job that is sent to a server with $I$ jobs is rejected.
Each incoming job is thus assigned to one out of $N$ queues as in Algorithm~\ref{alg}.
Upon each job arrival, we assume that the actions of sampling $d$ servers and assigning the job to some queue are instantaneous and occur at the same time.



Let $(Q^N(t),M^N(t))=(Q_k^N(t),M_k^N(t))_{k=1}^N\in \{0,1,\ldots,I\}_+^{2N}$ be the system state at time~$t\in\mathbb{R}_+$:~$Q_k^N(t)$ represents the number of jobs in queue $k$ at time $t$ and $M_k^N(t)$ represents the last observation collected from server~$k$ by the controller at time~$t$.
To avoid unnecessary technical complication in our proofs and since the observation associated to server~$k$ is no less than the actual number of jobs in~$k$ after sampling $k$ for the first time, for the initial condition we assume that $Q_k^N(0)\le M_k^N(0)$ for all $k$.

It is convenient to represent the system state by $X^N(t)=(X_{i,j}^N(t):0\le i\le j<\infty)$ where
\begin{equation}
X_{i,j}^N(t) \bydef \frac{1}{N}\sum_{k=1}^N \mathbf{1}_{\{ Q_k^N(t)=i, M_k^N(t)=j\}}
\end{equation}
denotes the proportion of $(i,j)$-servers at time $t$.
It is clear that $X^N(t)$ is still a Markov chain with values in some finite set $\mathcal{S}_N$ that is a subset of
\begin{equation}
\mathcal{S}\bydef \Big\{ \left(x_{i,j}\in\mathbb{R}_+ :0\le i\le j\le I\right): \sum_{i=0}^I \sum_{j=i}^I x_{i,j} =1 \Big\}.
\end{equation}
The transitions and rates of the Markov chain $X^N(t)$ that are due to server departures are easy to write because they have no impact on memory: for $x\in\mathcal{S}_N$, the transition $x\mapsto x - \frac{e_{i+1,j}}{N} + \frac{e_{i,j}}{N}$ occurs with rate~$N\,x_{i+1,j}$
where $e_{i,j}\bydef \left(\delta_{i,i'}\,\delta_{j,j'} \in\{0,1\} : 0\le i'\le j'\le I \right)$  and $\delta_{a,b}$ is the Kronecker delta.
%
On the other hand, the transitions and rates of  $X^N(t)$ that are due to job arrivals
are quite involved and they are omitted. However, in Section~\ref{sec:constructionX} we will show how to construct the sample paths of~$X^N(t)$.

\subsection{Fluid model.}
\label{eq:cus78dc}

For any $x\in\mathcal{S}$,
let
$$
x_{i,\cdot} \bydef \sum_{j= i}^I x_{i,j}
\quad \mbox{and} \quad 
x_{\cdot,j} \bydef \sum_{i=0}^j x_{i,j}
$$

The next definition introduces the \emph{fluid model} for the dynamics of $SQ(d,N)$.

\begin{definition}
\label{def1}
A function $x(t):\mathbb{R}_+\to \mathcal{S}$ is said to be a \emph{fluid model} (or fluid solution) if
the following conditions are satisfied:
\begin{enumerate}
 \item $x(t)$ is absolutely continuous, and
 \item $\frac{d x_{i,j}(t)}{dt}=b_{i,j}(x(t))$ almost everywhere, for every $i\ge 0$ and $j\ge i$,
\end{enumerate}
where
$b(x)\bydef (b_{i,j}(x):0\le i\le j\le I)$ is given by
\begin{align}
\label{eq:b00}
 b_{0,0}(x) =  &   \lambda d (x_{0,\cdot}-x_{0,0}) - \lambda  +  \mathcal{R}_0(x) \\
%
\nonumber
b_{i,j}(x) =  & x_{i+1,j}  -  \mathbf{1}_{\{i>0\}} x_{i,j} - \lambda d x_{i,j}
 - \mathcal{R}_{j-1}(x) \frac{x_{i,j}}{x_{\cdot,j}}  \,  \mathbf{1}_{\{ x_{\cdot,j} > 0 \}}\\
\label{eq:bij}
& + \mathbf{1}_{\{i>0\}} \mathcal{R}_{j-2}(x) \frac{x_{i-1,j-1}}{x_{\cdot,j-1}} \, \mathbf{1}_{\{ x_{\cdot,j-1} > 0 \}}
 + \mathbf{1}_{\{j=I,i>0\}} \mathcal{R}_{I-1}(x) \frac{x_{i-1,I}}{x_{\cdot,I}} \, \mathbf{1}_{\{x_{\cdot,I}>0\}},
 \qquad \forall i,j : i<j \\
 \label{eq:b11}
b_{1,1}(x) = &   -x_{1,1} + \lambda d (x_{1,\cdot}-x_{1,1})
 +  \lambda  -  \mathcal{R}_{0}(x)
  -  \mathcal{R}_{0}(x) \frac{x_{1,1}}{x_{\cdot,1}}  \, \mathbf{1}_{\{x_{\cdot,1}>0\}}   
  - \mathcal{G}_1 (x) \\
%
%
%
\nonumber
b_{i,i}(x) = &   -x_{i,i} + \lambda d (x_{i,\cdot}-x_{i,i})
  -  \mathcal{R}_{i-1}(x) \frac{x_{i,i}}{x_{\cdot,i}}  \, \mathbf{1}_{\{x_{\cdot,i}>0\}} 
  +  \mathcal{R}_{i-2}(x) \frac{x_{i-1,i-1}}{x_{\cdot,i-1}} \, \mathbf{1}_{\{x_{\cdot,i-1}>0\}}  \\
%
\label{eq:bii}
& + \mathcal{G}_{i-1} (x) - \mathcal{G}_{i} (x),
\qquad\qquad  \forall i=2,\ldots,I-1 \\
%
%
%
\label{eq:bII}
b_{I,I}(x) = &   -x_{I,I} 
  +  \mathcal{R}_{I-2}(x) \frac{x_{I-1,I-1}}{x_{\cdot,I-1}}  \, \mathbf{1}_{\{x_{\cdot,I-1}>0\}} 
%
+ \mathcal{G}_{I-1} (x)
+ \mathcal{R}_{I-1}(x) \frac{x_{I-1,I}}{x_{\cdot,I}}   \,\mathbf{1}_{\{x_{\cdot,I}>0\}}
\end{align}
with
\begin{equation}
\label{eq:structure_R}
%
%
%
\mathcal{R}_j(x) =  0 \vee \lambda  \left( 1 - d \sum_{i=0}^j  (j+1-i) \,x_{i,\cdot}    \right) \mathbf{1}_{\{\sum_{i=0}^j x_{\cdot,i}=0\}}
%
%
\end{equation}
\begin{equation}
\label{eq:Gi}
\mathcal{G}_j (x)
=  \lambda d   \,
\mathbf{1}_{\left\{ \sum_{i=0}^j x_{\cdot,i}=0,\, d \sum_{i=0}^j  (j+1-i) x_{i,\cdot} \le 1 \right\}}
\sum_{i=0}^j x_{i,\cdot}
\end{equation}
and $a\vee b \bydef \max\{a,b\}$.
\end{definition}

The discontinuous function $b$ will be referred to as \emph{drift}, and to some extent it may be interpreted as the conditional expected change from state $x$ of the Markov chain $X^N(t)$, though this may only be true when $x_{0,0}>0$, where $\mathcal{R}_{j}(x)=0$ for all $j$ and the formulas above become linear admitting a very intuitive explanation.


Let us provide some intuition for the drift expressions in Definition~\ref{def1}, and let us start with coordinates~(0,0).
At the moment of each arrival at the load balancer, the states of~$d$ servers are sampled and~$k$ idle servers that the load balancer has not yet spotted are sampled with probability ${ d\choose k } (x_{0,\cdot}-x_{0,0})^k ( 1 - x_{0,\cdot}+x_{0,0})^{d-k}$.
Since $  \sum_{k=1}^d k { d\choose k } (x_{0,\cdot}-x_{0,0})^k ( 1 - x_{0,\cdot}+x_{0,0})^{d-k} = d (x_{0,\cdot}-x_{0,0})$ and arrivals occur with rate~$\lambda$, the average rate in which $(0,j)$-servers, $j\ge 1$, are discovered is $\lambda d (x_{0,\cdot}-x_{0,0})$ and the rationale behind the first term in~\eqref{eq:b00} is justified.
The dynamics that remain to specify are the ones related to the effective job assignments, that is where singularities can happen.
%
In order to build a fluid model `consistent' with the finite stochastic system $X^N(t)$, 
one should take into account the fluctuations of order $1/N$ that appear when $X_{0,0}^N(t)=0$. These bring discontinuities in the drift.
%
Let $z_j  \bydef \sum_{i\ge j} x_{i,\cdot}$ and $R_j^N(t)\bydef \sum_{i=0}^j X_{\cdot,i}^N(t)$.
We notice that $\mathcal{R}_{j}(x)/\lambda$, where $\mathcal{R}_{j}(x)$ is defined in \eqref{eq:structure_R},
will be interpreted as the proportion of time where the process $(R_j^N(t))_{[t,t+\epsilon]}$ tends to stay on zero
with the load balancer sampling $(\cdot,j')$-servers only, for all $j'>j$, in the limit where $N\to\infty$ first and then $\epsilon\downarrow 0$; this will be formalized in Section~\ref{sec:fluid_any_coordinates}.
Thus, the term $\lambda -  \mathcal{R}_0(x)$ represents the rate in which jobs are assigned to (0,0)-servers, which become (1,1)-servers as soon as they receive a job.
%
This explains the drift expression in~\eqref{eq:b00}.
The particular structure of $\mathcal{R}_j(x)$ given in \eqref{eq:structure_R} will be the outcome of the stochastic analysis that will be developed in Section~\ref{sec:conn}.

Let us provide some intuition also for the drift expression on coordinates (1,1) (see \eqref{eq:b11}), as
it brings some additional interpretation that also applies on general coordinates.
The first term says that departures from $(1,1)$-servers occur with rate $x_{1,1}$ and the second one says that new $(1,1)$-servers are discovered with rate $\lambda d (x_{1,\cdot} -x_{1,1} )$.
This can be easily justified as done above for the first summation term of $b_{0,0}$.
Then, we notice that the $\lambda  -  \mathcal{R}_0(x)$ term has been already interpreted above and thus
the dynamics that remain to specify are the ones related to job assignments at (1,1)-servers.
According to SQ$(d,N)$, if the load balancer knows no $(0,0)$-server
then it randomizes over the set of $(\cdot,1)$-servers, and thus within this scenario $x_{1,1}$ should decrease with rate proportional to $\frac{x_{1,1}}{x_{\cdot,1}}$.
This is indeed the case if $x_{\cdot,1}>0$.
Thus, $\frac{x_{1,1}}{x_{\cdot,1}} \mathcal{R}_0(x)$ is the rate in which jobs are assigned to (1,1)-servers when $x_{\cdot,1}>0$.
It remains to model the rate in which jobs are assigned to (1,1)-servers when $x_{\cdot,1}=0$.
Since we aim at building a deterministic model `consistent' with the stochastic one, to model the rate of job assignments to (1,1)-servers when $x_{\cdot,1}=0$ one should take into account the fluctuations of order $1/N$ that appear when $R_{1}^N(t)=0$. 
The term
$\mathcal{G}_1 (x)$ given in \eqref{eq:Gi}
is indeed such rate, and again it will be the outcome of the stochastic analysis developed in Section~\ref{sec:conn}.



The following proposition will be proven in Section~\ref{sec:conn}.

\begin{proposition}
\label{thexistence}
Fluid solutions exist.
\end{proposition}

\subsection{Connecting the Markov and the fluid models.}
\label{simulations}

Our first result is the following connection between the stochastic and the fluid models.


\begin{theorem}
\label{th0}
Assume that $X^N(0)\to x^0\in\mathcal{S}$ almost surely.
With probability one, any limit point of the stochastic process $(X^{N}(t))_{t\in[0,T]}$
satisfies the conditions that define a fluid solution. 
\end{theorem}
In view of this result, proven in Section~\ref{sec:conn}, a fluid solution may be interpreted as an accurate approximation of the time-dependent dynamics of the finite stochastic system $X^N(t)$, provided that $N$ is sufficiently large. 

Given $x\in\mathcal{S}$, let us define the functions
\begin{equation}
\nonumber
\mathcal{L}_S(x) \bydef \sum_{i=1}^I i x_{i,\cdot}
\qquad
\mathcal{L}_M(x) \bydef \sum_{j=1}^I j x_{\cdot,j}
\end{equation}
We notice that
$\mathcal{L}_S(X^{N}(t))$
represents the number of jobs in the system at time~$t$ scaled by~$N$
and that
$\mathcal{L}_M(X^{N}(t))$
represents the number of jobs scaled by~$N$ the load balancer \emph{believes} are in the system at time~$t$.
Since the system is symmetric with respect to the servers, the function $\mathcal{L}_S(X^{N}(t))$ is also interpreted as the average number of jobs at time $t$ in each queue.

It is clear that
$\mathcal{L}_M(x) - \mathcal{L}_S(x) =  \sum_{i\ge 0}\sum_{j\ge i}  (j-i) x_{i,j} \ge 0,$
which is to be expected because ($\cdot,j$)-servers can not contain more than~$j$ jobs by definition.

The following corollary of Theorem~\ref{th0} is immediate.
\begin{corollary}
\label{cor1}
Let $x(t)$ be a fluid solution.
Assume that $(x(t))_{t\in[0,T]}$ is a limit point of $(X^{N}(t))_{t\in[0,T]}$ with probability one.
Then, 
$(\mathcal{L}_S(x(t)))_{t\in[0,T]}$
and 
$(\mathcal{L}_M(x(t)))_{t\in[0,T]}$
are limit points of
$(\mathcal{L}_S(X^{N}(t)))_{t\in[0,T]}$
and
$(\mathcal{L}_M(X^{N}(t)))_{t\in[0,T]}$, respectively, with probability one.
\end{corollary}



\begin{figure}
\label{fig:simulations}
\centering
\includegraphics[width=\columnwidth]{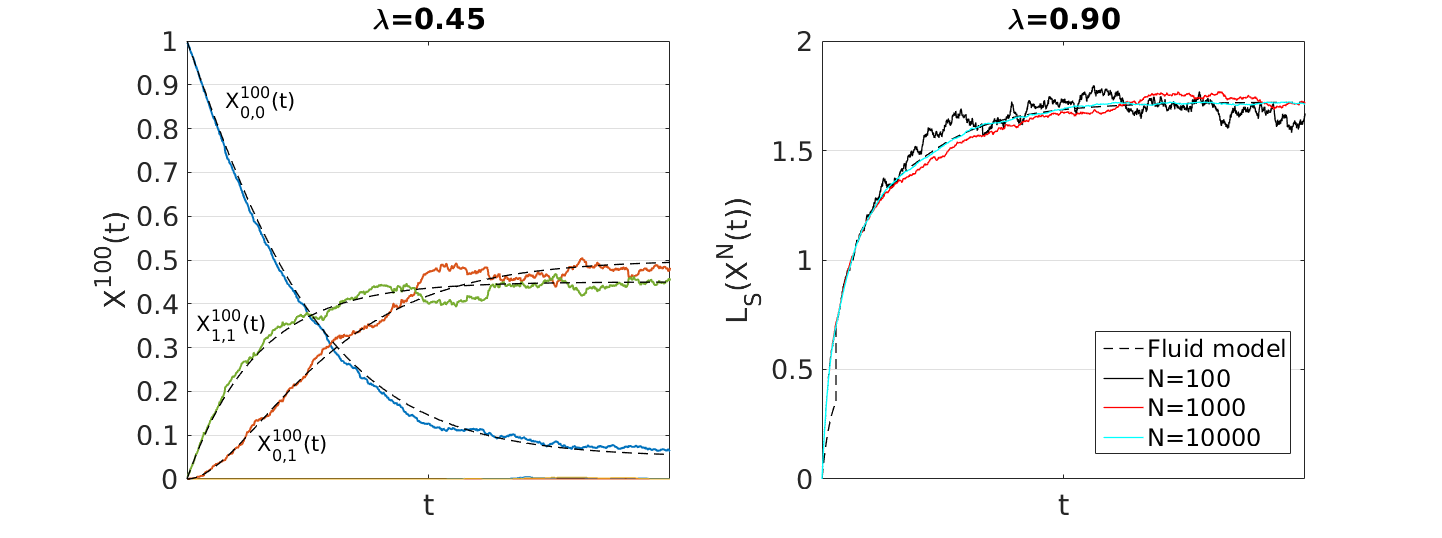}
\caption{Numerical convergence of the stochastic model $X^N(t)$ (continuous lines) to the fluid model $x(t)$ (dashed line).}
\end{figure}

We complement Theorem~\ref{th0} and Corollary~\ref{cor1} presenting some numerical simulations to support the claim that the fluid model provides a remarkably accurate approximation of the sample paths of $X^N(t)$ even when~$N$ is finite and relatively small.
%
Assuming $d=2$, Figure~\ref{fig:simulations} plots the time dependent dynamics of $X^{N}(t)$ and $x(t)$.
At time zero, we have chosen $X^{N}(0)$ and $x(0)$ such that $X_{0,0}^{N}(0) = x_{0,0}(0)=1$, which means that all servers are idle and the load balancer is aware of it.
Each curve on these plots is an average over ten simulations.
The fluid (stochastic) model is always represented by dashed (continuous) lines.
In the picture on the left ($\lambda=0.45$), we set $N=100$ and notice that the fluid model already captures in an accurate manner the dynamics of $X^N(t)$, which turn out to be concentrated more and more on just three components: namely (0,0), (0,1) and (1,1).
Matter of fact $X_{0,0}^{N}(t)+X_{0,1}^{N}(t)+X_{1,1}^{N}(t)$ gets closer and closer to 1 when both $N$ and $t$ increase.
In the picture on the right ($\lambda=0.9$), dynamics are distributed on several components and for convenience
we have plotted $\mathcal{L}_S(X^{N}(t))$ and its fluid model counterpart $\mathcal{L}_S(x(t))$.
We notice that $\mathcal{L}_S(x(t))$ almost overlaps the trajectory of $\mathcal{L}_S(X^{N}(t))$ already when $N=1000$.
This size is in agreement with the magnitude of modern distributed computing such as web-server farms or data-centers, as they are often composed of (tenths of) thousands of servers.

\section{Main results}
\label{sec:matching}


In this section we focus on fluid solutions and investigate optimality and stability properties.
First, we are interested in fixed points.
\begin{definition}
We say that a fluid solution $x(t)$ is a \emph{fixed point} if $b(x(t))=0$ for all $t$.
\end{definition}

When fluid solutions are fixed points, we drop the dependency on $t$.

Let us define  $j^\star$ as in \eqref{eq:jstar}
and for simplicity let us assume that $I>j^\star$.

The next result, proven in Section~\ref{proof_th1}, establishes the existence and uniqueness of a fixed point and says that its mass is concentrated only on coordinates of the form $(i,j^\star)$ and $(i,j^\star+1)$.

\begin{theorem}[Existence and Uniqueness of Fixed Points]
\label{th1}
There exists a unique fixed point, say $x^\star$.
It is such that $x_{\cdot, j^\star}^\star+x_{\cdot, j^\star+1}^\star=1$ and
\begin{subequations}
\begin{align}
\label{eq:x_0j}
\lambda d \, x_{0,j^\star}^\star   &= (1+\lambda d  ) (1-\lambda) - \frac{1}{(1+ \lambda d)^{j^\star}} \\
\label{eq:x_0j1}
x_{0,j^\star}^\star + x_{0,j^\star+1}^\star &=1-\lambda.
%
\end{align}
\end{subequations}
\end{theorem}

In the fixed point, our first remark is that queue lengths are bounded, by $j^\star+1$.
%
%
As we show in our proof, an explicit expression for $x^\star$ seems to be difficult to obtain, though it can be easily computed when $\lambda$ and $d$ are fixed numerically. 
In fact, in Section~\ref{proof_th1} we provide an explicit expression for $x_{i,j}^\star$ when $(i,j)\neq (j^\star,j^\star)$
as a function of $x_{j^\star,j^\star}^\star$,
and identify $x_{j^\star,j^\star}^\star$ by means of a polynomial equation of degree $j^\star+1$ (see \eqref{eq:xjj}).

A case of particular interest is when $j^\star=0$, which given \eqref{eq:jstar} occurs if and only if $\lambda<1-1/d$, where we have the following remark.

\begin{remark}[Asymptotic Optimality]
\label{remark1}
If $\lambda<1-1/d$,
then Theorem~\ref{th1} implies that $x_{0,0}^\star=1-\lambda-1/d$, $x_{0,1}^\star=1/d$, $x_{1,1}^\star=\lambda$ and $x_{i,j}^\star=0$ on the remaining coordinates.
Thus, provided that dynamics converge to $x^\star$,
we have shown that a load balancer implementing SQ($d,N$) is always aware of the fact that some servers are idle when $ N\to\infty$ and $t$ is sufficiently large because  $x_{0,0}^\star>0$.
In this scenario, the load balancer can certainly assign each incoming job to one of such idle servers, and the job itself would incur zero delay.
This is in fact the ideal situation for any arriving job and in this sense we say that SQ$(d,N)$ is \emph{asymptotically optimal}.
\end{remark}



%

The next proposition provides further insights on the system performance at the fixed point $x^\star$.

\begin{proposition}
\label{prop:LS_LM}
Let $x^\star$ as in Theorem~\ref{th1}. Then,
\begin{equation}
\label{skcj878dc}
\mathcal{L}_M(x^\star)=\mathcal{L}_S(x^\star) + \frac{1}{d}
\end{equation}
and
\begin{equation}
\label{skcj878dcsvdf}
j^\star  -\frac{1}{d} \le \mathcal{L}_S(x^\star) \le  j^\star  -\frac{1}{d} + 1.
\end{equation}
\end{proposition}

Proposition~\ref{prop:LS_LM}, proven in Section~\ref{section:LS_LM}, provides simple bounds on the average number of jobs in each queue.
It also says that there is a fluid mass equal to ${1}/{d}$ that the load balancer will never spot.
In other words, the samplings performed by the load balancer at each arrival will correctly build the true state of the system up to an (absolute) error of~$1/d$.


In Remark~\ref{remark1}, we discussed the asymptotic optimality of SQ($d,N$) postulating some form of stability for fluid solutions when $t\to\infty$.
The next result shows that fluid solutions are indeed globally stable and that convergence to $x^\star$ occurs exponentially fast, provided that $\lambda < 1-1/d$.


\begin{theorem}[Global Stability]
\label{th2}
Let $x(t)$ be a fluid solution.
If $\lambda <1-1/d$, then
there exist $\alpha>0$ and $\beta>0$ independent of~$t$ such that
\begin{equation}
\label{eqsjauis}
\|x(t)-x^\star\|\le \alpha e^{-\beta t},\quad \forall t
\end{equation}
where $\|\cdot\|$ is the Euclidean norm.
\end{theorem}
The proof of Theorem~\ref{th2} is given in Section~\ref{proof_th2} and is based on the following `Lyapunov-type' argument.
When $x_{0,0}(t)=0$, we first show that $\dot{\mathcal{L}}_S(x(t)) \le \lambda -1+\frac{1}{d}$, which implies
that $\mathcal{L}_S(x(t))$ decreases with derivative bounded away from zero.
However, since $\mathcal{L}_S(x(t))\ge 0$, $x_{0,0}(t)$ must necessarily increase in finite time, and when it does we show that $x(t)$ is uniquely determined by the unique solution of 
a linear ODE system of the form $\dot{x}=A(x-x^\star)$.
At this point, \eqref{eqsjauis} follows by standard results of ODE theory.
%
When $\lambda \ge 1-1/d$, a generalization of this argument 
is complicated by the involved structure of~$x^\star$ and the fact that $\mathcal{L}_S(x(t))$ is in general not monotone.
However, we conjecture that $x^\star$ remains globally stable.
This is also confirmed by the numerical simulations shown in Section~\ref{simulations}.

\section{Connection between the fluid and the Markov models}
\label{sec:conn}


We now prove that the sequence of stochastic processes  $\{(X^{N}(t))_{t\in[0,T]}\}_{N=d}^\infty$ converges almost surely, as $N\to\infty$, to a fluid solution, for any $T>0$.
This proves Proposition~\ref{thexistence} and Theorem~\ref{th0}.

Our proof is based on three steps.
First, we construct the sample paths of the process $X^{N}(t)$ on each pair of coordinates. This is achieved using a common coupling technique that defines the processes $(X^{N}(t))_{t\in[0,T]}$ for all $N\in\mathbb{Z}_+$ on a single probability space and in terms of a finite number of ``fundamental processes''.
Then, we show that limit trajectories exist and are Lipschitz continuous with probability one. This is done by using standard arguments, e.g., Gamarnik et al.~\cite{Gamarnik2016}, Tsitsiklis and Xu~\cite{tsitsiklis2012}, and Bramson~\cite{Bramson1998}.
Finally, we prove that any such limit trajectory must be a fluid solution, which is the main difficulty.
This last step is based on technical arguments that are specific to the stochastic model under investigation.

\subsection{Probability space and coupled construction of sample paths.}
\label{sec:constructionX}

We construct a probability space where the stochastic processes $\{(X^{N}(t))_{t\in[0,T]}\}_{N\ge d}$ are coupled.
All the processes of interest will be a function of the following fundamental processes, all of them independent of each other:

\begin{itemize}
\item $\mathcal{N}_{\lambda} (t)$, the Poisson processes of job arrivals, with rate~$\lambda$,
defined on $(\Omega_A, \mathcal{A}_A, \mathbb{P}_A)$;

\item $\mathcal{N}_{1} (t)$, the Poisson processes of potential job departures, with rate 1,
defined on $(\Omega_D, \mathcal{A}_D, \mathbb{P}_D)$;

\item $V_n^p$ for all $p=1,\ldots,d$, $(W_n)_n$, $(U_n)_n$,
where
the random variables $V_n^p$, $W_n$ and $U_n$, for all $n$ and $p$, are all independent and uniformly distributed over the interval $(0,1]$. These are \emph{selection} processes: 
$(V_n^p)_n$ will select the servers to sample at each arrival (see Line~5 of Algorithm~\ref{alg}),
$(W_n)_n$ will be used to randomize among the servers having the lowest observations (see Line~8 of Algorithm~\ref{alg})
and 
$(U_n)_n$ will select the server that fires a departure.
These $2+d$ processes are defined on $(\Omega_S, \mathcal{A}_S, \mathbb{P}_S)$;

%

\item $(X^N(0))_N$, the process of the initial conditions, where each random variable $X^N(0)$ takes values in $\mathcal{S}$,
defined on $(\Omega_0, \mathcal{A}_0, \mathbb{P}_0)$.

\end{itemize}


Each process $\{(X^{N}(t))_{t\in[0,T]}\}$, with $N\ge d$, can be constructed on
$(\Omega,\mathcal{A},\mathbb{P}) =  (\Omega_A \times\Omega_D \times \Omega_S \times\Omega_0, \mathcal{A}_A\times \mathcal{A}_D\times \mathcal{A}_S\times \mathcal{A}_0, \mathbb{P}_A \times \mathbb{P}_D \times\mathbb{P}_S \times\mathbb{P}_0)$
by using that $\mathcal{N}_{\lambda} (N t)=_{st}\mathcal{N}_{\lambda N} (t)$, where $=_{st}$ denotes equality in distribution.
This equality ensures that the Poisson process with rate $\lambda N$, which represents the arrival process associated to the $N$-th system, is coupled with the fundamental Poisson process~$\mathcal{N}_{\lambda} (t)$.
Since $\mathcal{N}_{1} (N t)=_{st}\mathcal{N}_{N} (t)$, this coupling is also used for the processes of potential job departures.



Now, let $t_n^{N,\lambda}$ and $t_n^{N,1}$ be the times of the $n$-th jump of the Poisson processes $\mathcal{N}_\lambda(N t)$ and $\mathcal{N}_1(N t)$, respectively.
Let also $X_{i,j}^N(t^-)\bydef \lim_{s \uparrow t } X_{i,j}^N(s)$ and $X_{i,\cdot}^N (t) \bydef \sum_{j\ge i} X_{i,j}^N (t)$.
In view of the coupling discussed above, we can construct
$X_{0,0}^N(t)$ as follows
\begin{subequations}
\label{X00Nt_aa}
\begin{align} 
\label{X00Nt}
 X_{0,0}^N(t) =   X_{0,0}^N(0)&
+ \frac{1}{N} \sum_{n=1}^{\mathcal{N}_\lambda(N t)} \sum_{p=1}^d  \mathbb{I}_{ ( X_{0,0}^N(t_n^{N,\lambda-}), X_{0,\cdot}^N(t_n^{N,\lambda-})   ] }^{(V_n^p)} \\
\label{X00Nt_1}
& + \frac{1}{N} \sum_{n=1}^{\mathcal{N}_\lambda(N t)}
\left( \mathbf{1}_{ \{ X_{0,0}^N(t_n^{N,\lambda-})=0 \} } \prod\limits_{p=1}^d  \mathbb{I}_{ (   X_{0,\cdot}^N(t_n^{N,\lambda}-), 1   ] }^{(V_n^p)}   - 1  \right).
\end{align}
\end{subequations}
In the above expression,
the term \eqref{X00Nt} corresponds to the action of sampling $d$ servers
and 
the term \eqref{X00Nt_1} corresponds to the action of assigning each job to a server.
At the arrival of the $n$-th job, $t_n^{N,\lambda}$,\
the proportion of (0,0)-servers
increases by $k/N$ if $k$
$(0,j)$-servers, for any $j>0$, are sampled, which justifies the term in~\eqref{X00Nt},
and
decreases by $1/N$ except when such proportion is zero immediately before $t_n^{N,\lambda}$ and no idle server is sampled at time $t_n^{N,\lambda}$, which justifies the term in~\eqref{X00Nt_1}.
%

Using the random variables $W_n$ and $U_n$, an expression similar to \eqref{X00Nt_aa} can be written for $X_{i,j}^N(t)$ when $(i,j)\in\{0,1,\ldots,I\}^2$.
Towards this purpose,
let us define
\begin{equation}
\label{eq:R_S_Z}
R_{i}^N(t) \bydef  \sum_{j=0}^{i}  X_{\cdot,j}^N(t), \quad
S_{i,j}^N(t) \bydef  \sum_{i'=0}^{i-1}  X_{i',\cdot}^N(t) + \sum_{j'\ge i}^j  X_{i,j'}^N(t), \quad
 Z_{i}^N(t) \bydef  \sum_{k\ge i}  X_{k,\cdot}^N(t),
\end{equation}
which respectively represent
\emph{i)}  the proportion of servers that the controller believes have at most~$i$ jobs,
\emph{ii)} the proportion of servers with at most~$i-1$ jobs, or $i$ jobs but with observation less than or equal to~$j$
and \emph{iii)} the proportion of servers with at least~$i$ jobs, 
and
\begin{equation}
\nonumber
%
%
M_{i,j,n}^N \bydef \frac{1}{N}\sum_{p=1}^d  \mathbb{I}_{ ( S_{i,j-1}^N(t_n^{N,\lambda-}),  S_{i,j}^N(t_n^{N,\lambda-})   ] }^{(V_n^p)}, \quad
\underline{M}_{j,n}^N \bydef \sum_{i=0}^j M_{i,j,n}^N,
 \quad
\overline{M}_{j,n}^N \bydef \sum_{j'\ge j} M_{j',j,n}^N.
\end{equation}
We notice that $M_{i,j,n}^N$, $\underline{M}_{j,n}^N$ and $\overline{M}_{j,n}^N$ are the scaled-by-$N$ numbers of $(i,j)$-, $(\cdot,j)$- and $(j,\cdot)$- servers sampled immediately before time $t_n^{N,\lambda}$, respectively.
Furthermore, 
let also 
\begin{equation}
\label{def:F_ijn}
F_{i,j,n}^N\bydef 
\mathbb{I}_{ \left( \sum_{k=0}^{i-1} X_{k,j}^N(t_n^{N,\lambda-}) -M_{k,j,n}^N,  \sum_{k=0}^{i} X_{k,j}^N(t_n^{N,\lambda-}) - M_{k,j,n}^N  \right] } 
^{(W_n (X_{\cdot,j}^N(t_n^{N,\lambda-}) + \overline{M}_{j,n}^N - \underline{M}_{j,n}^N )  )}
%
\end{equation}
if $i<j$, and
\begin{equation}
\label{def:F_iin}
F_{j,j,n}^N\bydef 
\mathbb{I}_{
\left( 
\sum_{k=0}^{j-1} X_{k,j}^N(t_n^{N,\lambda-}) -M_{k,j,n}^N,
X_{\cdot,j}^N(t_n^{N,\lambda-}) + \overline{M}_{j,n}^N - \underline{M}_{j,n}^N
\right]
}
^{(W_n (X_{\cdot,j}^N(t_n^{N,\lambda-}) + \overline{M}_{j,n}^N - \underline{M}_{j,n}^N )  )}
%
\end{equation}
if $j\ge 1$.
For all $i\le j$, the random variable $F_{i,j,n}^N$ will be used to handle the randomness in Line~8 of Algorithm~\ref{alg} and thus perform a job assignment to a $(i,j)$-server, which needs to be chosen in the set of $(\cdot,j)$-servers.
Specifically, we will use $F_{i,j,n}^N$, with $i\le j$, in the scenario where $R_{i-1}^N(t_n^{N,\lambda-})=0$ and $M_{i',j',n}^N =0$ for all $i'< i$ and $j'\ge i'$, that is the case where
the load balancer memory contains no observation less than $j$
and no server containing less than $j$ jobs is sampled immediately before~$t_n^{N,\lambda}$.
In this case, according to SQ$(d,N)$, the $n$-th job must be routed to a random $(\cdot,j)$-server, provided that
such a server exists.
This randomness is captured by the uniform random variable $W_n$ and we notice that 
$N (X_{\cdot,j}^N(t_n^{N,\lambda-}) + \overline{M}_{j,n}^N - \underline{M}_{j,n}^N )$ is the number of 
$(\cdot,j)$-servers, or equivalently the occurrences of $j$ in the memory of the load balancer, at the arrival of the $n$-th job  and after having performed the associated sampling of the states of $d$ random servers.
Within these conditions, the job arriving at time $t_n^{N,\lambda}$ is routed to an $(i,j)$-server if and only if $F_{i,j,n}^N=1$.

Provided that $i<j$, the following formula constructs the process $X^N(t)$ on coordinates $(i,j)$
\begin{subequations}
\label{eq:shucyu}
\begin{align}
%
\label{eq:shucyu1}
 X_{i,j}^N(t) =  X_{i,j}^N(0) & + 
\frac{1}{N} \sum_{n=1}^{\mathcal{N}_1(N t)} \mathbb{I}_{ ( S_{i+1,j}^N(t_n^{N,1-})-X_{i+1,j}^N(t_n^{N,1-}),  S_{i+1,j}^N(t_n^{N,1-})   ] }^{(U_n)} \\
\label{eq:shucyu1bis}
& - \frac{ \mathbf{1}_{\{i>0\}} }{N} \sum_{n=1}^{\mathcal{N}_1(N t)}  \mathbb{I}_{ ( S_{i,j}^N(t_n^{N,1-})-X_{i,j}^N(t_n^{N,1-}),  S_{i,j}^N(t_n^{N,1-})   ] }^{(U_n)}\\
& 
\label{eq:shucyu2}
- \frac{1}{N} \sum_{n=1}^{\mathcal{N}_\lambda(N t)}
\sum_{p=1}^d  \mathbb{I}_{ ( S_{i,j}^N(t_n^{N,\lambda-})-X_{i,j}^N(t_n^{N,1-}),  S_{i,j}^N(t_n^{N,\lambda-})   ] }^{(V_n^p)}
\\
\label{eq:shucyu3}
& - \frac{1}{N} \sum_{n=1}^{\mathcal{N}_\lambda(N t)} \mathbf{1}_{\{R_{j-1}^N(t_n^{N,\lambda-})=0\}} F_{i,j,n}^N \prod\limits_{p=1}^d \mathbb{I}_{ ( 1-Z_j^N(t_n^{N,\lambda-}), 1 ] } ^{(V_n^p)} \\
\label{eq:shucyu4}
& + \frac{\mathbf{1}_{\{i>0\}}}{N} \sum_{n=1}^{\mathcal{N}_\lambda(N t)} \mathbf{1}_{\{R_{j-2}^N(t_n^{N,\lambda-})=0\}} F_{i-1,j-1,n}^N \prod\limits_{p=1}^d \mathbb{I}_{ ( 1-Z_{j-1}^N(t_n^{N,\lambda-}), 1 ] }^{(V_n^p)}\\
\label{eq:shucyu5}
& + \frac{\mathbf{1}_{\{j=I,i>0\}}}{N} \sum_{n=1}^{\mathcal{N}_\lambda(N t)} \mathbf{1}_{\{R_{I-1}^N(t_n^{N,\lambda-})=0\}} F_{i-1,I,n}^N \prod\limits_{p=1}^d \mathbb{I}_{ ( 1-Z_{I}^N(t_n^{N,\lambda-}), 1 ] }^{(V_n^p)}.
\end{align}
\end{subequations}
The summations in~\eqref{eq:shucyu1} and~\eqref{eq:shucyu1bis} refer, respectively, to job departures from $(i+1,j)$- and $(i,j)$-servers,
the summation in~\eqref{eq:shucyu2} refers to the case where $k$
$(i,j)$-servers are sampled (as soon as $k$ of them are sampled, they become $(i,i)$-servers, and thus $X_{i,j}^N$ decreases by $k/N$), and
the summations in~\eqref{eq:shucyu3} and~\eqref{eq:shucyu4} refer to the case where a job is assigned to an $(i,j)$-server and to an $(i-1,j-1)$-server, respectively.
We notice that
a job can be assigned at time $t_n^{N,\lambda}$ to an $(i,j)$-server only if 
the memory contains no server with observation less than~$j-1$ immediately before $t_n^{N,\lambda}$ (i.e., $R_{j-1}^N(t_n^{N,\lambda-})=0$)
and no $(i',j')$-server, for some $i'<i$ and for any $j'$, has been sampled at time $t_n^{N,\lambda}$.
Summation \eqref{eq:shucyu5} covers the boundary case where $j=I$ and has the same intuition of term~\eqref{eq:shucyu5}.

Similarly, 
when $i=j\ge 1$, we have
\begin{subequations}
\label{eq:shucyuascs}
\begin{align} 
\label{eq:scuas78cas1}
 X_{i,i}^N(t) = X_{i,i}^N(0) & - \frac{1}{N} \sum_{n=1}^{\mathcal{N}_1(N t)}  \mathbb{I}_{ ( S_{i-1,I}^N(t_n^{N,1-}),  S_{i,i}^N(t_n^{N,1-})   ] }^{(U_n)}\\
\label{eq:scuas78cas2}
& +  \frac{1}{N} \sum_{n=1}^{\mathcal{N}_\lambda(N t)}
\sum_{p=1}^d  \mathbb{I}_{ ( S_{i,i}^N(t_n^{N,\lambda-}), S_{i,I}^N(t_n^{N,\lambda-})   ] }^{(V_n^p)}
\\
%
%
%
\label{eq:scuas78cas3}
& + \frac{\mathbf{1}_{\{i=1\}}}{N} \sum_{n=1}^{\mathcal{N}_\lambda(N t)}  \mathbf{1}_{\{X_{0,0}^N(t_n^{N,\lambda-}) + \sum_{j=1}^I M_{0,j,n}^N>0\}}\\
\label{eq:scuas78cas3bis}
& + \frac{\mathbf{1}_{\{i>1\}}}{N} \sum_{n=1}^{\mathcal{N}_\lambda(N t)}  \mathbf{1}_{ \left\{R_{i-2}^N(t_n^{N,\lambda-}) = 0 \right\} } F_{i-1,i-1,n}^N  \prod\limits_{p=1}^d \mathbb{I}_{ ( 1-Z_{i-1}^N(t_n^{N,\lambda-}), 1 ] }^{(V_n^p)} \\
%
%
%
\label{eq:scuas78cas4}
&  - \frac{\mathbf{1}_{\{i<I\}}}{N} \sum_{n=1}^{\mathcal{N}_\lambda(N t)}  \mathbf{1}_{ \left\{R_{i-1}^N(t_n^{N,\lambda-}) = 0 \right\} } F_{i,i,n}^N 
\prod\limits_{p=1}^d \mathbb{I}_{ ( 1-Z_{i}^N(t_n^{N,\lambda-}), 1 ] }^{(V_n^p)} \\
\label{eq:scuas78cas5}
& + \frac{\mathbf{1}_{\{i=I\}}}{N} \sum_{n=1}^{\mathcal{N}_\lambda(N t)} \mathbf{1}_{\{R_{I-1}^N(t_n^{N,\lambda-})=0\}} F_{I-1,I,n}^N \prod\limits_{p=1}^d \mathbb{I}_{ ( 1-Z_{I}^N(t_n^{N,\lambda-}), 1 ] }^{(V_n^p)}.
\end{align}
\end{subequations}
The summation in~\eqref{eq:scuas78cas1} refers to job departures from $(i,i)$-servers,
the summation in~\eqref{eq:scuas78cas2} refers to the sampling of $k$ different $(i,j)$-servers, which become $(i,i)$-servers immediately after sampling.
Finally, the summations in~\eqref{eq:scuas78cas3} and~\eqref{eq:scuas78cas4} refer to jobs assignments and have the same intuition of~\eqref{eq:shucyu3} and~\eqref{eq:shucyu4}.


\subsection{Limit trajectories are Lipschitz.}

With respect to a set of sample paths $\omega$ having probability one, we show that any subsequence of the sequence $\{X^N(\omega,t)\}_N$ contains a further subsequence $\{X^{N_k}(\omega,t)\}_{k}$ that converges to some Lipschitz continuous function $x$.
This proves tightness of sample paths. 

First, let us introduce the following formulas for quick reference.
These can be proven in a straightforward manner using the strong law of the large numbers and the functional strong law of large numbers for the Poisson process.

\begin{lemma}
\label{GC_F}
Let $T>0$ and $a,b\in [0,1]^d$ such that $a_k\le b_k$ for all $k=1,\ldots,d$.
There exists $\mathcal{C}\subseteq \Omega$ such that $\mathbb{P}(\mathcal{C})=1$
such that
\begin{equation}
\nonumber
\lim_{N\to\infty} \sup_{t\in[0,T]} | \tfrac{1}{N} \mathcal{N}_\lambda(N t,\omega) -\lambda t| =0,
\qquad
\lim_{N\to\infty} \sup_{t\in[0,T]} | \tfrac{1}{N} \mathcal{N}_1(N t,\omega) - t| =0,
\end{equation}
\begin{equation}
\nonumber
\lim_{N\to\infty} \frac{1}{N} \sum_{n=1}^N 
\sum_{p=1}^d  \mathbb{I}_{ ( a_k,b_k  ] }^{(V_n^p(\omega))}
= \sum_{p=1}^d  b_k-a_k,
\qquad
\lim_{N\to\infty} \frac{1}{N} \sum_{n=1}^N 
\prod_{p=1}^d  \mathbb{I}_{ ( a_p,b_p  ] }^{(V_n^p(\omega))}
= \prod_{p=1}^d  (b_p-a_p)
\end{equation}
for all $\omega\in\mathcal{C}$.
\end{lemma}

In the following, we will work on the set $\mathcal{C}$ introduced in previous lemma and
we will also often use that
\begin{equation}
\nonumber
\lim_{N\to\infty} \frac{1}{N} \sum_{n=1}^{\mathcal{N}_\lambda(Nt,\omega)} 
\prod_{p=1}^d  \mathbb{I}_{ ( a_k,b_k  ] }^{(V_n^p(\omega))}
= \lambda t \prod_{p=1}^d  (b_k-a_k)
\end{equation}
by the renewal theorem.

%
%
%
%

Let $x^0\in[0,1]$, sequences $R_n\downarrow 0$ and $\gamma_n\downarrow 0$, and a constant $L>0$ be given.
For $n\ge 1$, let also
\begin{subequations}
\begin{align*}
E_N(R_N, \gamma_N,L,x^0) \bydef \big\{ & x\in D[0,T] :  |x(0)-x^0|\le R_N,  \, |x(a)-x(b)| \le L|a-b|+\gamma_N, \,\forall a,b\in[0,T]  \big\}
\end{align*}
\end{subequations}
and
\begin{subequations}
\begin{align*}
E_c(L,x^0) \bydef \big\{ & x\in D[0,T] : x(0)=x^0,  \, |x(a)-x(b)| \le L|a-b|, \,\forall a,b\in[0,T]  \big\}.
\end{align*}
\end{subequations}

The next lemma says that the sample paths along any coordinates $(i,j)$ is approximately Lipschitz continuous.
The proof is omitted because follows exactly the same standard arguments used in Lemma~5.2 of Gamarnik et al. \cite{tsitsiklis2017}, which basically use the fact that the jumps of the Markov chain of interest are of the order of $1/N$ and that the evolution of such Markov chain on a given pair of coordinates only depends on the evolution of such Markov chain on a finite number of other coordinates.

\begin{lemma}
\label{LM:tightness}
Fix $T>0$, $\omega\in\mathcal{C}$, and some $x^0\in\mathcal{S}$.
Suppose that $\|X^N(\omega,0)-x^0\| \le \tilde{R}_N$,
for some sequence $\tilde{R}_N\downarrow 0$. Then, there exists sequences $R_N\downarrow 0 $ and $\gamma_N\downarrow 0$ such that
\begin{equation}
\nonumber
X_{i,j}^N(\omega,\cdot)\in E_N(R_N, \gamma_N,L,x^0),\quad \forall (i,j)\in\mathbb{Z}^+: i\le j, \forall N
\end{equation}
where  $L=\lambda d+1$.
\end{lemma}

The next proposition says that the sample paths along any coordinates $(i,j)$ are sufficiently close to a Lipschitz continuous function. 
The proof is omitted because follows exactly the same arguments used in the proof of Proposition~11 in Tsitsiklis and Xu \cite{tsitsiklis2012}: it uses Lemma~\ref{LM:tightness} 
and topological properties of the space $E_c(L,x^0)$, i.e., sequential compactness (by the Arzel\`{a}-Ascoli theorem) and closedness.


\begin{proposition}
\label{PR:tightness}
Fix $T>0$, $\omega\in\mathcal{C}$, and some $x^0\in\mathcal{S}$.
Suppose that $\|X^N(\omega,0)-x^0\| \le \tilde{R}_N$,
for some sequence $\tilde{R}_N\downarrow 0$. Then,
every subsequence of $\{X^N(\omega,\cdot )\}_{N=1}^\infty$ contains a further subsequence 
$\{X^{N_k}(\omega,\cdot )\}_{k=1}^\infty$ such that
\begin{equation}
\nonumber
\lim_{k\to\infty} \sup_{t\in[0,T]} |X_{i,j}^{N_k} (\omega, t) - x_{i,j}(t)| =0, \quad \forall i,j\ge i
\end{equation}
where $x_{i,j}\in E_c(1+\lambda d,x^0)$ for all $i,j\ge i$.
\end{proposition}

Since Lipschitz continuity implies absolute continuity, we have thus obtained that limit points of $X^N(t)$ are absolutely continuous, and
it remains to show that the partial derivatives of $x(t)$ are given by the expressions in Definition~\ref{def1}.


\subsection{Limit trajectories are fluid solutions.}
\label{ascascs}

To conclude the proof of Theorem~\ref{th0}, it remains to show that any limit point is a fluid solution, i.e., it satisfies the conditions given in Definition~\ref{def1}. 
This is the main technical difficulty.

Fix $\omega\in\mathcal{C}$ and let $\{X^{N_k}(\omega,t)\}_{k=1}^\infty$ be a subsequence that converges to $\overline{x}$, i.e.
\begin{equation}
\label{seb_conv}
\lim_{k\to\infty} \sup_{t\in [0,T]} \| X^{N_k}(\omega, t) - \overline{x}(t)\| = 0.
\end{equation}
Since $\overline{x}_{i,j}$ must be Lipschitz continuous for all $i$ and $j$ by Proposition~\ref{PR:tightness}, it is also absolutely continuous and thus it remains to show that 
\begin{equation}
\label{eqhcus}
\dot{\overline{x}}_{i,j}(t) =\lim_{\epsilon\to 0} \frac{1}{\epsilon} \lim_{k\to\infty} X_{i,j}^{N_k}(t+\epsilon) - X_{i,j}^{N_k}(t) = b_{i,j}(\overline{x}(t)),
\end{equation}
whenever $\overline{x}_{i,j}(\cdot)$ is differentiable.
This will be done in the following subsections. Now, we introduce the following technical lemmas.

\begin{lemma}
\label{lemma:X_close_x}
Fix $\omega\in\mathcal{C}$, $\epsilon>0$ and let \eqref{seb_conv} hold.
Then, 
for all $i, j$ and $ t$,
$$|X_{i,j}^{N_k}(u)-\overline{x}_{i,j}(t)| \le 2 L \epsilon,\quad \forall u \in [t,t+\epsilon]$$ 
for all $k$ sufficiently large, where $L=\lambda d+1$.
\end{lemma}
\emph{Proof:}
By Lemma~\ref{LM:tightness}, there exists a sequence $\gamma_{N_k}\downarrow 0$ such that 
$X_{i,j}^{N_k}(\omega,u) \in [ \overline{x}_{i,j}(t) -\epsilon L - \gamma_{N_k}, \overline{x}_{i,j}(t) +\epsilon L + \gamma_{N_k} ]$, 
for all $u\in[t,t+\epsilon]$.
Thus, for all $k$ sufficiently large,
$X_{i,j}^{N_k}(\omega,u) \in [ \overline{x}_{i,j}(t) -2\epsilon L, \overline{x}_{i,j}(t) +2\epsilon L  ]$,
 for all $u\in[t,t+\epsilon]$, as desired.
\hfill$\Box$
\indent

As a corollary of Lemma~\ref{lemma:X_close_x}, we obtain
\begin{equation}
\label{X_close_x}
|S_{i,j}^{N_k}(u)-s_{i,j}(\overline{x}(t))| \le C \epsilon,\quad \forall u \in [t,t+\epsilon]
\end{equation}
for all $k$ sufficiently large, where $C\bydef 2L(I+1)^2$.

\begin{remark}
In the following, we will work on any fixed trajectory $\omega\in\mathcal{C}$ but we will write $X^{N_k}(t)$, instead of $X^{N_k}(\omega,t)$, for simplicity of notation.
\end{remark}

\begin{lemma}
\label{lemma_ascasji1}
Fix $\omega\in\mathcal{C}$ and let \eqref{seb_conv} hold.
Then, 
\begin{subequations}
\label{lemma_ascasji12}
\begin{align}
\label{lemma_ascasji12AAA}
 \lim_{\epsilon\to 0}\lim_{k\to\infty} \frac{1}{\epsilon N_k} \sum_{n=\mathcal{N}_\lambda(N_k t)+1}^{\mathcal{N}_\lambda(N_k (t+\epsilon))}
\sum_{p=1}^d  \mathbb{I}_{ ( S_{i,j}^{N_k}(t_n^{N_k,\lambda-}) - X_{i,j}^{N_k}(t_n^{N_k,\lambda-}), S_{i,j}^{N_k}(t_n^{N_k,\lambda-})   ] }^{(V_n^p)}
 &\,=  \lambda d  \overline{x}_{i,j}(t) \\
 \label{lemcj_hcyty}
 \lim_{\epsilon\to 0}\lim_{k\to\infty} \frac{1}{\epsilon N_k} \sum_{n=\mathcal{N}_1(N_k t)+1}^{\mathcal{N}_1(N_k (t+\epsilon))}
\mathbb{I}_{ ( S_{i,j}^{N_k}(t_n^{N_k,\lambda-}) - X_{i,j}^{N_k}(t_n^{N_k,\lambda-}), S_{i,j}^{N_k}(t_n^{N_k,\lambda-})   ] }^{(U_n)}
 &\,=  \overline{x}_{i,j}(t).
\end{align}
\end{subequations}
\end{lemma}
\emph{Proof:}
Given in the Appendix.
\hfill$\Box$
\indent

\subsubsection{Fluid solution on coordinates (0,0).}
\label{sec:fluid00}

The next lemma explicits the derivative of $\overline{x}_{0,0}(t)$ when $\overline{x}_{0,0}(t)>0$. It also implies that $\overline{x}_{0,0}(\cdot)$ is differentiable when strictly positive.

\begin{lemma}
\label{lemma:x00_whenx00positive}
Fix $\omega\in\mathcal{C}$, let \eqref{seb_conv} hold and assume $\overline{x}_{0,0}(t)>0$.
Then,
\begin{equation}
\label{asjc89as}
\dot{\overline{x}}_{0,0}(t)
=-  \lambda   + d \lambda  (\overline{x}_{0,\cdot}(t) - \overline{x}_{0,0}(t)).
\end{equation}
\end{lemma}
\emph{Proof:}
Choose $\epsilon>0$ small enough such that $\overline{x}_{0,0}(t) -2(I+1)\epsilon L>0$ where $L=\lambda d+1$.
Such  $\epsilon$ exists because $\overline{x}_{0,0}(t)>0$ by hypothesis.
Since $t_n^{N_k,\lambda-}\in (t,t+\epsilon]$ when $n \in \{ \mathcal{N}_\lambda({N_k} t)+1, \ldots, \mathcal{N}_\lambda({N_k} (t+\epsilon))\}$, 
Lemma~\ref{lemma:X_close_x} implies that
for all~$k$ sufficiently large
we must have
$$
\mathbf{1}_{ \{ X_{0,0}^{N_k}(t_n^{N_k,\lambda}-)>0 \} } = 1, \quad \forall n \in \{ \mathcal{N}_\lambda({N_k} t)+1, \ldots, \mathcal{N}_\lambda({N_k} (t+\epsilon))\}.
$$
Thus, using~\eqref{X00Nt_aa}, we obtain
\begin{equation}
\nonumber
 \lim_{k\to\infty} X_{0,0}^{N_k}(t+\epsilon) - X_{0,0}^{N_k}(t) =  \lim_{k\to\infty}
 \frac{1}{N_k} \sum_{n=\mathcal{N}_\lambda({N_k} t)+1}^{\mathcal{N}_\lambda({N_k} (t+\epsilon))} 
  -1 + \sum_{p=1}^d  \mathbb{I}_{ ( X_{0,0}^{N_k}(t_n^{{N_k},\lambda-}), X_{0,\cdot}^{N_k}(t_n^{{N_k},\lambda-})   ] }^{(V_n^p)}.
\end{equation}
A direct application of Lemma~\ref{lemma_ascasji1} concludes the proof.
\hfill$\Box$
\indent

%
%
%
%
%

The next two lemmas give properties on the boundary where $\overline{x}_{0,0}(t)=0$.

\begin{lemma}
\label{lemma:x00_not_differentiable}
Fix $\omega\in\mathcal{C}$, let \eqref{seb_conv} hold and assume
$\overline{x}_{0,0}(t)=0$ and $d \overline{x}_{0,\cdot}(t)> 1$. Then,
$t$ is not a point of differentiability.
\end{lemma}
\emph{Proof:}
First of all, we notice that
\begin{equation}
\label{asjuyyrhfv}
\lim_{\epsilon\downarrow 0} \frac{1}{\epsilon} \lim_{k\to\infty} X_{0,0}^{N_k}(t+\epsilon) - X_{0,0}^{N_k}(t)
%
\ge \\
-  \lambda   + d \lambda  (\overline{x}_{0,\cdot}(t) - \overline{x}_{0,0}(t)).
\end{equation}
This holds true because
\begin{equation}
\label{scja78sc}
 X_{0,0}^{N_k}(t+\epsilon) - X_{0,0}^{N_k}(t)
%
\ge 
  \frac{1}{N_k} \sum_{n=\mathcal{N}_\lambda(N_k t)+1}^{\mathcal{N}_\lambda(N_k (t+\epsilon))} 
\sum_{p=1}^d  \mathbb{I}_{ ( X_{0,0}^{N_k}(t_n^{N_k,\lambda-}), X_{0,\cdot}^{N_k}(t_n^{N_k,\lambda-})   ] }^{(V_n^p)} - 1,
\end{equation}
which is obvious given~\eqref{X00Nt_aa}, and because the RHS of~\eqref{scja78sc}, once divided by $\epsilon$, converges to $-  \lambda   + d \lambda  (\overline{x}_{0,\cdot}(t) - \overline{x}_{0,0}(t))$, by Lemmas~\ref{GC_F} and~\ref{lemma_ascasji1}, in the limit where $k\to\infty$ first and $\epsilon\downarrow 0$.

Now, assume by contradiction that $t$ is a point of differentiability.
In this case, the limit in the LHS of \eqref{asjuyyrhfv} exists and must be equal to $\dot{\overline{x}}_{0,0}(t)$.
Furthermore, if  $\overline{x}_{0,\cdot}(t) > 1/d$, the RHS of~\eqref{asjuyyrhfv} is strictly positive and thus $\dot{\overline{x}}_{0,0}(t)$ must be strictly positive as well.
On the other hand, it is not possible to have $\dot{\overline{x}}_{0,0}(t)>0$ and $\overline{x}_{0,0}(t)=0$ because the function $\overline{x}_{0,0}$ is always non-negative. This contradicts that~$t$ is a point of differentiability of $\overline{x}_{0,0}(\cdot)$.
\hfill$\Box$
\indent

The next lemma says that the limit trajectory $\overline{x}_{0,0}$ remains on zero in a right neighborhood of~$t$, provided that $\overline{x}_{0,0}(t)=0$ and $0\le\overline{x}_{0,\cdot}(t)< 1/d$.

\begin{lemma}
\label{lemma:x00_stay_at_zero}
Fix $\omega\in\mathcal{C}$, let \eqref{seb_conv} hold and assume
$\overline{x}_{0,0}(t)=0$ and $d \overline{x}_{0,\cdot}(t)< 1$. Then,
\begin{equation}
\label{ascj876}
 \exists \delta>0 : \overline{x}_{0,0}(t')=0,\,\, \forall t'\in[t,t+\delta].
\end{equation}
\end{lemma}
\emph{Proof:}
Assume that \eqref{ascj876} is false.
Then, there exists a sequence $t_n\downarrow t$ such that $t_n>t_{n+1}>t$ and
\begin{equation}
\nonumber
\overline{x}_{0,0}(t_n)>0
\mbox{ and } 
\dot{\overline{x}}_{0,0}(t_n)>0
\end{equation}
for all $n$.
By Lemma~\ref{lemma:x00_whenx00positive}, we have
$\dot{\overline{x}}_{0,0}(t_n) = -  \lambda   + d \lambda  (\overline{x}_{0,\cdot}(t_n) - \overline{x}_{0,0}(t_n))$
and thus
$\overline{x}_{0,\cdot}(t_n)-\overline{x}_{0,0}(t_n) > \frac{1}{d}$, for all $n$,
and by continuity
\begin{equation}
\inf_n \overline{x}_{0,\cdot}(t_n)-\overline{x}_{0,0}(t_n) \ge \frac{1}{d}.
\end{equation}
Thus, we get
\begin{equation}
\nonumber
\overline{x}_{0,\cdot}(t)-\overline{x}_{0,0}(t)
=     \lim_{n\to\infty} \overline{x}_{0,\cdot}(t_n)-\overline{x}_{0,0}(t_n) 
\ge   \inf_n \overline{x}_{0,\cdot}(t_n)-\overline{x}_{0,0}(t_n) 
\ge   \frac{1}{d}.
\end{equation}
This contradicts the hypothesis.
\hfill$\Box$
\indent

Summarizing,
\begin{itemize}
 \item  when $ \overline{x}_{0,0}(t) >0$, we have proven that $\dot{\overline{x}}_{00}(t) =b_{0,0}(\overline{x}(t))$;
 \item  when $ \overline{x}_{0,\cdot}(t) < 1/d$ and $\overline{x}_{0,0}(t) = 0$, we have proven that $\overline{x}_{0,0}(t)$ remains zero on a right neighborhood, and thus if $t$ is a point of differentiability, then $0 = \dot{\overline{x}}_{00}(t) = b_{0,0}(\overline{x}(t))$;
 \item  when $ \overline{x}_{0,\cdot}(t) > 1/d$ and $\overline{x}_{0,0}(t) = 0$, we have proven that $t$ is not a point of differentiability;
 \item  when $ \overline{x}_{0,\cdot}(t) = 1/d$ and $\overline{x}_{0,0}(t) = 0$, either $t$ is not a point of differentiability or it is. In the latter case, we must have $\dot{\overline{x}}_{00}(t) =0 $ because $\overline{x}_{00}$ is a non-negative function and since also $b_{0,0}(\overline{x}(t))=0$, we have indeed $\dot{\overline{x}}_{00}(t) =b_{0,0}(\overline{x}(t))$ as desired.
\end{itemize}
Thus, $\dot{\overline{x}}_{0,0}(t) = b_{0,0}(\overline{x}(t))$ almost everywhere.

\subsubsection{Fluid solution on arbitrary coordinates.}
\label{sec:fluid_any_coordinates}

We now prove that $\dot{\overline{x}}_{i,j}(t) = b_{i,j}(\overline{x}(t))$ almost everywhere with respect to arbitrary coordinates $(i,j)$.
This requires a more in-depth analysis of the stochastic process $X^N(t)$.

Let 
\begin{equation}
\label{eq:R_j_x}
R_j(t)
\bydef 
\lim_{\epsilon \downarrow 0} \lim_{k\to\infty} \frac{1}{\epsilon N_k }
\sum_{n=\mathcal{N}_\lambda(N_k t)+1}^{\mathcal{N}_\lambda(N_k (t+\epsilon))} 
\mathbf{1}_{ \left\{R_{j}^{N_k}(t_n^{N_k,\lambda-}) = 0 \right\} }
\prod\limits_{p=1}^d \mathbb{I}_{ ( 1- Z_{j+1}^{N_k}(t_n^{N_k,\lambda-}), 1 ] }^{(V_n^p)}
,
\end{equation}
which is interpreted as the proportion of time where the process $R_{j}^{N_k}(t_n^{N_k,\lambda-})$ remains on zero in the interval $[t,t+\epsilon]$
while the load balancer keeps sampling only $(j',\cdot)$-servers for all $j'>j$
in the $\lim_{\epsilon \downarrow 0} \lim_{k\to\infty}$ limit.
In the following, we show that $R_j(t)=\mathcal{R}_j(\overline{x}(t))$, where $\mathcal{R}_j$ is given in Definition~\ref{def1}.

The structure of $R_0$ is easily obtained as a corollary of the analysis developed in previous section.

\begin{lemma}
\label{prop_R0}
Fix $\omega\in\mathcal{C}$, let \eqref{seb_conv} hold, and assume that $\overline{x}(t)$ is differentiable.
Then, $R_0(t)$ exists and is given by
\begin{equation}
\label{eq:R_0}
R_0(t) = 
\lambda (1 -  d \,\overline{x}_{0,\cdot}(t)) \, \mathbf{1}_{\{\overline{x}_{0,0}(t)=0  \}} \mathbf{1}_{\{d \,\overline{x}_{0,\cdot}(t)) < 1 \}}.
\end{equation}
\end{lemma}
\emph{Proof:}
First, we notice that if ${\overline{x}}_{0,0}(t)>0$, then necessarily $R_0(t)=0$.
In fact, if for any $j$, $\sum_{i=0}^j \overline{x}_{\cdot,i}(t)>0$, then we can find $\epsilon>0$ such that $\sum_{i=0}^j \overline{x}_{\cdot,i}(t)-2L(I+1)^2\epsilon>0$.
Since $t_n^{N_k,\lambda-}\in (t,t+\epsilon]$ when $n \in \{ \mathcal{N}_\lambda({N_k} t)+1, \ldots, \mathcal{N}_\lambda({N_k} (t+\epsilon))\}$, 
Lemma~\ref{lemma:X_close_x} implies that
for all~$k$ sufficiently large
we must have $ R_{j}^{N_k}(t_n^{N_k,\lambda-})>0$, and therefore $R_{j}(t)=0$ in this case.

Thus, assume that $\overline{x}_{0,0}(t)=0$.
In this case, since $t$ is a point of differentiability, we necessarily have $\dot{\overline{x}}_{0,0}(t)=0$ and, by Lemma~\ref{lemma:x00_not_differentiable}, necessarily $d \,\overline{x}_{0,\cdot}(t) \le 1$. This gives the indicator functions in~\eqref{eq:R_0}.
Furthermore, recalling the structure of $X_{0,0}^N(t)$ given in~\eqref{X00Nt_aa}, 
we have
\begin{align*}
&0= \dot{\overline{x}}_{0,0}(t)\\
&=  \lim_{\epsilon \downarrow 0} \frac{1}{\epsilon} \lim_{k\to\infty} X_{0,0}^{N_k}(t+\epsilon)-X_{0,0}^{N_k}(t)  \\
&= 
\lim_{\epsilon \downarrow 0} \frac{1}{\epsilon} \lim_{k\to\infty} 
\frac{1}{N_k} \sum_{n=1+\mathcal{N}_\lambda(N_k t)}^{\mathcal{N}_\lambda(N_k (t+\epsilon))} 
 \mathbf{1}_{ \{ X_{0,0}^{N_k}(t_n^{{N_k},\lambda-})=0 \} } \prod\limits_{p=1}^d  \mathbb{I}_{ (   X_{0,\cdot}^{N_k}(t_n^{{N_k},\lambda}-), 1   ] }^{(V_n^p)}   - 1 
+ \sum_{p=1}^d  \mathbb{I}_{ ( X_{0,0}^{N_k}(t_n^{N_k,\lambda-}), X_{0,\cdot}^{N_k}(t_n^{N_k,\lambda-})   ] }^{(V_n^p)}\\
&= 
\lambda d (\overline{x}_{0,\cdot}(t) - \overline{x}_{0,0}(t)) - \lambda + 
R_0(t).
\end{align*}
In the last equality we have used Lemma~\ref{lemma_ascasji1} and the definition of $R_0$.
This equation gives~\eqref{eq:R_0}.
\hfill$\Box$
\indent

The next lemma provides an expression for $R_j(t)$ for all $j$ and shows that $R_j(t)=\mathcal{R}_j(\overline{x}(t))$.
Our proof, given in the appendix, is based on Lemma~\ref{prop_R0}, which allows us to establish the existence and find the structure of $R_j$ in an iterative manner.

\begin{lemma}
\label{lemma:R}
Fix $\omega\in\mathcal{C}$,  let \eqref{seb_conv} hold and assume that $\overline{x}(t)$ is differentiable.
Then, for all $j$, $R_j(t)$ exists and is given by
\begin{equation}
\label{eq:sdjviuy7r1}
R_{j}(t) = 0 \vee \lambda  \left( 1 - d \sum_{i=0}^j  (j+1-i) x_{i,\cdot}(t)    \right) \mathbf{1}_{\{\sum_{i=0}^j x_{\cdot,i}(t)=0\}}.
\end{equation}
\end{lemma}

For any $i,j\ge i$, let us define
\begin{equation}
\Gamma_{i,j}^{\epsilon,k}(t) \bydef 
\frac{1}{\epsilon N_k}
\sum_{n=\mathcal{N}_\lambda(N_k t)+1}^{\mathcal{N}_\lambda(N_k (t+\epsilon))} 
\mathbf{1}_{ \left\{R_{j}^{N_k}(t_n^{N_k,\lambda-}) = 0 \right\} }
F_{i,j+1,n}^{N_k}
\prod_{p=1}^d \mathbb{I}_{(1 - Z_{j+1}^{N_k}(t_n^{N_k,\lambda -}) ,1]}^{(V_n^p)}
\end{equation}
and $\Gamma_{i,j}(t) \bydef  \lim_{\epsilon \downarrow 0}  \lim_{k\to\infty} \Gamma_{i,j}^{\epsilon,k}(t)$,
which is interpreted as the proportion of time where the process $R_{j}^{N_k}(t_n^{N_k,\lambda-})$ remains on zero in the interval $[t,t+\epsilon]$
while the load balancer
samples $(j',\cdot)$-servers only, for all $j'> j$,
and 
assigns jobs to $(i,j+1)$-servers only when the proportion of $(\cdot,j+1)$-servers vanishes
in the $\lim_{\epsilon \downarrow 0} \lim_{k\to\infty}$ limit.

The next lemma, proven in the appendix, gives an expression for $\Gamma_{i,j}(t)$
when $\overline{x}_{\cdot,j+1}(t)>0$
and will allow us to identify the limit behavior of terms  \eqref{eq:shucyu3}-\eqref{eq:shucyu5} and \eqref{eq:scuas78cas3bis}-\eqref{eq:scuas78cas5}.
\begin{lemma}
\label{lemma:iv8fd}
Fix $\omega\in\mathcal{C}$ and let \eqref{seb_conv} hold.
Assume that $\overline{x}(t)$ is differentiable and that $\overline{x}_{\cdot,j+1}(t)>0$. Then,
\begin{equation}
\label{asjc78dcckiujff}
\Gamma_{i,j}(t) = \frac{\overline{x}_{i,j+1}(t)}{\overline{x}_{\cdot,j+1}(t)} R_j(t).
\end{equation}
for all $i,j$ such that $ i\le j+1$.
\end{lemma}

With the lemmas above, we can identify the asymptotic behavior (in the $\lim_{\epsilon\downarrow 0} \lim_{k\to\infty}$ limit) of each summation appearing in the expressions 
of $X_{i,j}^{N_k}(t+\epsilon) - X_{i,j}^{N_k}(t)$
and 
$X_{i,i}^{N_k}(t+\epsilon) - X_{i,i}^{N_k}(t)$
that are obtained using \eqref{eq:shucyu} and \eqref{eq:shucyuascs}, respectively.

Let us first treat the case $i<j$.

Applying Lemma~\ref{lemma:iv8fd} in \eqref{eq:shucyu}, when $\overline{x}(t)$ is differentiable we obtain
\begin{subequations}
\begin{align}
\nonumber
\dot{\overline{x}}_{i,j}(t) = &\, \lim_{\epsilon\to 0} \frac{1}{ \epsilon} \lim_{k\to\infty}  X_{i,j}^{N_k}(t+\epsilon) - X_{i,j}^{N_k}(t) \\
\nonumber
=&\,  \overline{x}_{i+1,j}(t)  -  \mathbf{1}_{\{i>0\}} \overline{x}_{i,j}(t) - \lambda d \overline{x}_{i,j}(t)\\
\nonumber
 & - \frac{\overline{x}_{i,j}(t)}{\overline{x}_{\cdot,j}(t)} R_{j-1}(t)  \,  \mathbf{1}_{\{ \overline{x}_{\cdot,j}(t) > 0 \}}\\ 
\nonumber
& + \mathbf{1}_{\{i>0\}} \frac{\overline{x}_{i-1,j-1}(t)}{\overline{x}_{\cdot,j-1}(t)} R_{j-2}(t)\, \mathbf{1}_{\{ \overline{x}_{\cdot,j-1}(t)>0 \}}\\
%
\nonumber
& + \mathbf{1}_{\{j=I,i>0\}} \frac{\overline{x}_{i-1,I}(t)}{\overline{x}_{\cdot,I}(t)} R_{I-1}(t) \, \mathbf{1}_{\{\overline{x}_{\cdot,I}(t)>0\}}\\
\label{cskdjuyghhjh}
&+\lim_{\epsilon \downarrow 0}  \lim_{k\to\infty}  - \Gamma_{i,j-1}^{\epsilon,k}(t) \mathbf{1}_{\{\sum_{j'=0}^{j} \overline{x}_{\cdot,j'}(t) =0 \}} + \Gamma_{i-1,j-2}^{\epsilon,k}(t) \mathbf{1}_{\{i>0, \sum_{j'=0}^{j-1} \overline{x}_{\cdot,j'}(t) =0 \}}
%
\end{align}
\end{subequations}
where the first three terms follow by applying Lemma~\ref{lemma_ascasji1} to terms~\eqref{eq:shucyu1}, \eqref{eq:shucyu1bis} and \eqref{eq:shucyu2}.
%
Now, assume $i=0$ and $j>0$. 
Then, if $\sum_{j'=0}^{j}\overline{x}_{\cdot,j'}(t) =0$, then the first six terms of previous equation and the second summation term in \eqref{cskdjuyghhjh} are equal to zero,
and if in addition $t$ is a point of differentiability, then necessarily $\dot{\overline{x}}_{i,j}(t)=0$ (because $\overline{x}_{i,j}(t)=0$), which means that necessarily
$$
\lim_{\epsilon \downarrow 0}  \lim_{k\to\infty}  \Gamma_{i,j-1}^{\epsilon,k}(t) \mathbf{1}_{\{\sum_{j'=0}^{j} \overline{x}_{\cdot,j'}(t) =0 \}} 
$$
exists and is equal to zero.
Assume $i>0$ and $j>i$. 
If $\sum_{j'=0}^{j}\overline{x}_{\cdot,j'}(t) =0$, then the first six terms of previous equation again coincide with zero
and if in addition $t$ is a point of differentiability, then necessarily 
\begin{equation}
\label{sckj8uigb}
\lim_{\epsilon \downarrow 0}  \lim_{k\to\infty}  - \Gamma_{i,j-1}^{\epsilon,k}(t) \mathbf{1}_{\{\sum_{j'=0}^{j} \overline{x}_{\cdot,j'}(t) =0 \}} + \Gamma_{i-1,j-2}^{\epsilon,k}(t) \mathbf{1}_{\{i>0, \,\sum_{j'=0}^{j-1} \overline{x}_{\cdot,j'}(t) =0 \}}
\end{equation}
exists and is equal to zero.
Furthermore, if $\sum_{j'=0}^{j-1}\overline{x}_{\cdot,j'}(t) =0$ and $\overline{x}_{\cdot,j}(t) > 0$, then \eqref{sckj8uigb} still exists and is equal to zero as a consequence of the fact that we have inductively shown that $\lim_{\epsilon \downarrow 0}  \lim_{k\to\infty}  \Gamma_{i-1,j-2}^{\epsilon,k}(t) \mathbf{1}_{\{i>0, \overline{x}_{\cdot,j-1}(t) =0 \}} =0$.
Therefore, the limit in \eqref{cskdjuyghhjh} is always equal to zero.
We have thus shown that $\dot{\overline{x}}_{i,j}=b_{i,j}(\overline{x}_{i,j})$ when $i<j$.

The case $i=j>0$ is treated in a similar manner.
Let $G_i^{\epsilon,k}(t) \bydef \mathbf{1}_{ \{\sum_{j'=0}^{i} \overline{x}_{\cdot,j'}(t)=0 \}} \Gamma_{i,i-1}^{\epsilon,k}(t)$ 
and  $G_i \bydef G_i(t) \bydef \lim_{\epsilon\downarrow 0} \lim_{k\to\infty} G_i^{\epsilon,k}(t)$.
In the following, we show that $G_i(t)=\mathcal{G}_i(\overline{x}(t))$, where $\mathcal{G}_i$ is given in Definition~\ref{def1}.

Applying Lemma~\ref{lemma_ascasji1} to handle terms \eqref{eq:scuas78cas1} and \eqref{eq:scuas78cas2}, 
rewriting  term \eqref{eq:scuas78cas3} as
\begin{align*}
\mathbf{1}_{\{X_{0,0}^N(t_n^{N,\lambda-}) + \sum_{j=0}^I M_{0,j,n}^N>0\}}
= &1 - \mathbf{1}_{\{X_{0,0}^N(t_n^{N,\lambda-})=0\}}  \mathbf{1}_{\{\sum_{j=0}^I M_{0,j,n}^N=0\}} \\
= &1 - \mathbf{1}_{ \left\{X_{0,0}^{N_k}(t_n^{N_k,\lambda-}) = 0 \right\} }
\prod_{p=1}^d \mathbb{I}_{(1 - Z_{1}^{N_k}(t_n^{N_k,\lambda -}) ,1]}^{(V_n^p)},
\end{align*}
and applying Lemma~\ref{lemma:iv8fd} to handle terms \eqref{eq:scuas78cas3bis}, \eqref{eq:scuas78cas4} and \eqref{eq:scuas78cas5},
when $\overline{x}(t)$ is differentiable we obtain
\begin{subequations}
\label{eq:s78casc}
\begin{align}
\dot{\overline{x}}_{i,i}(t) =  &  \lim_{\epsilon\to 0} \frac{1}{\epsilon} \lim_{k\to\infty} X_{i,i}^{N_k}(t+\epsilon) - X_{i,i}^{N_k}(t) \\
\label{eq:s78cascgdt3_A0}
 = &   -\overline{x}_{i,i}(t) + \lambda d (\overline{x}_{i,\cdot}(t)-\overline{x}_{i,i}(t))\\
\label{eq:s78cascgdt3_A}
& + \mathbf{1}_{\{i=1\}}  (\lambda -  R_0(t)z_1^d )  \\
\label{eq:s78cascgdt3_B}
& + \mathbf{1}_{\{i>1\}} 
R_{i-2}(t) \frac{\overline{x}_{i-1,i-1}(t)}{\overline{x}_{\cdot,i-1}(t)} \, \mathbf{1}_{\{\overline{x}_{\cdot,i-1}(t)>0\}} 
\\
\label{eq:s78cascgdt3_C}
&  - 
\mathbf{1}_{\{i<I\}} 
R_{i-1}(t) \frac{\overline{x}_{i,i}(t)}{\overline{x}_{\cdot,i}(t)}  \, \mathbf{1}_{\{\overline{x}_{\cdot,i}(t)>0\}}
\\
\label{eq:s78cascgdt3_D}
& + \mathbf{1}_{\{i=I\}} R_{I-1}(t) \frac{\overline{x}_{I-1,I}(t)}{\overline{x}_{\cdot,I}(t)}  \, \mathbf{1}_{\{\overline{x}_{\cdot,I}(t)>0\}}\\
\label{eq:s78cascgdt3_E}
&+ \lim_{\epsilon\downarrow 0} \lim_{k\to\infty}   - G_i^{\epsilon,k}(t) \mathbf{1}_{ \{i<I \} }  + G_{i-1}^{\epsilon,k}(t) \mathbf{1}_{ \{i>1 \} } .
\end{align}
\end{subequations}
Now, assume that $i=1$.
If $t$ is a point of differentiability and $\overline{x}_{0,0}(t)+\overline{x}_{\cdot,1}(t)=0$, then we must have $\dot{\overline{x}}_{1,1}(t)=0$, and thus
necessarily
\begin{align*}
G_1(t) =  \left( \lambda d  \overline{x}_{1,\cdot}(t)  + \lambda -  R_0(t)
\right)  \mathbf{1}_{\{ \overline{x}_{0,0}(t)+\overline{x}_{\cdot,1}(t)=0 \}}
= \lambda d  \left( \overline{x}_{1,\cdot}(t)  + \overline{x}_{0,\cdot}(t)
\right)  \mathbf{1}_{\{ \overline{x}_{0,0}(t)+\overline{x}_{\cdot,1}(t)=0 \}}.
\end{align*}
In the last equality we have used Lemma~\ref{prop_R0} and that $d\overline{x}_{0,\cdot}(t)\le 1$, which holds true because $t$ is a point of differentiability (Lemmas~\ref{lemma:x00_not_differentiable} and~\ref{lemma:x00_stay_at_zero}).
When $i=2,\ldots,I-1$,
if $t$ is a point of differentiability and $\sum_{j'=0}^{i} \overline{x}_{\cdot,j'}(t)=0$,
then necessarily $\dot{\overline{x}}_{i,i}(t)=0$ and proceeding in an iterative manner, we obtain
\begin{align*}
G_i(t) & = \left(  G_{i-1}(t) + \lambda d \overline{x}_{i,\cdot}(t) \right) \mathbf{1}_{\{ \sum_{i'=0}^i \overline{x}_{\cdot,i'}(t)=0 \}} 
=  \lambda d  \mathbf{1}_{\{ \sum_{i'=0}^i \overline{x}_{\cdot,i'}(t)=0 \}} \sum_{i'=0}^i \overline{x}_{i',\cdot}(t).
\end{align*}


The following two lemmas are a generalization of Lemmas~\ref{lemma:x00_not_differentiable} and~\ref{lemma:x00_stay_at_zero} and show under which conditions $\overline{x}(t)$ is differentiable. The proofs use the same arguments in those lemmas and therefore they are omitted.

\begin{lemma}
\label{lemma:xdotj_not_differentiable}
Fix $\omega\in\mathcal{C}$, let \eqref{seb_conv} hold and assume
$\sum_{j'=0}^j\overline{x}_{\cdot,j'}(t)=0$ and $d \sum_{i=0}^j  (j+1-i) \overline{x}_{i,\cdot}(t) > 1$. Then,
$t$ is not a point of differentiability.
\end{lemma}

\begin{lemma}
\label{lemma:xdotj_stay_at_zero}
Fix $\omega\in\mathcal{C}$, let \eqref{seb_conv} hold and assume
$\sum_{j'=0}^j\overline{x}_{\cdot,j'}(t)=0$ and $d \sum_{i=0}^j  (j+1-i) \overline{x}_{i,\cdot}(t) < 1$. Then,
\begin{equation}
 \exists \delta>0 : \sum_{j'=0}^j\overline{x}_{\cdot,j'}(t')=0,\,\, \forall t'\in[t,t+\delta].
\end{equation}
\end{lemma}

Now, we notice that the expressions of $R_{j}(t) $ and $G_i(t)$ obtained so far assumed that~$t$ was a point of differentiability.
However, if $\sum_{j'=0}^j\overline{x}_{\cdot,j'}(t)=0$, previous lemmas say that this can only be true if $d \sum_{i=0}^j  (j+1-i) x_{i,\cdot}(t) \le 1$.
Thus, those expressions make sense only in that case.
This does not change the structure of $R_{j}(t)$ obtained in \eqref{eq:sdjviuy7r1} because
\begin{align*}
 R_{j}(t) \times \mathbf{1}_{\left\{ d \sum_{i=0}^j  (j+1-i) \overline{x}_{i,\cdot}(t) \le 1 \right\}} = 
 R_{j}(t)
%
= \mathcal{R}_j(\overline{x}(t)).
\end{align*}
but on the other hand we must have
$G_i(t)=\mathcal{G}_i(\overline{x}(t))$, where $\mathcal{G}_i$ is defined in \eqref{eq:Gi}.
We have thus shown that $\dot{\overline{x}}_{i,i}(t)=b_{i,i}(x)$.


\section{Proofs of Theorems~\ref{th1} and~\ref{th2}}
\label{proofs_fluid}

Let us introduce the intervals
\begin{equation}
\nonumber
\mathcal{I}_n\bydef \left[ \lambda_{n}^*, \lambda_{n+1}^* \right),\quad \forall n\ge 0
\end{equation}
where $\lambda_0^*=0$ and $\lambda_n^*$, for $n\ge 1$, is the unique root in $(0,1]$ of the polynomial equation
\begin{equation}
\nonumber
(1-z)(z d + 1)^n = 1.
\end{equation}
The first values of $\lambda_n^*$ are  $\lambda_1^*=1-\frac{1}{d}$ and $\lambda_2^*=\frac{1}{2}-\frac{1}{d} +\sqrt{\frac{1}{4}+\frac{1}{d}}$.
We notice that $j^\star$, defined in \eqref{eq:jstar}, is the unique integer such that $\lambda \in \mathcal{I}_{j^\star}$.
In fact, $\lambda \in \mathcal{I}_{n}$ if and only if 
\begin{align}
\nonumber
n=  - \frac{ \log (1-\lambda_n^*)}{\log (\lambda_n^* d +1)} \le - \frac{ \log (1-\lambda)}{\log (\lambda d +1)}
\end{align}
and
\begin{align}
\nonumber
n+1=  - \frac{ \log (1-\lambda_{n+1}^*)}{\log (\lambda_{n+1}^* d +1)} > - \frac{ \log (1-\lambda)}{\log (\lambda d +1)},
\end{align}
which thus implies $n = \left\lfloor  - \frac{ \log (1-\lambda)}{\log (\lambda d +1)} \right\rfloor =j^\star$.

Let $x(t)$ be a fluid solution.
Since the $x_{i,j}(t)$'s are absolutely continuous, both $\mathcal{L}_S(x(t))$ and $\mathcal{L}_M(x(t))$ are absolutely continuous as well, and thus almost everywhere differentiable.
When~$t$ is a point of differentiability, it is clear that 
\begin{equation}
\label{eq:Lyapunov_function_derivative}
\dot{\mathcal{L}}_S(x(t)) =  \sum_{i=1}^I i \,b_{i,\cdot}(x(t)),\qquad
\dot{\mathcal{L}}_M(x(t)) =  \sum_{j=1}^I j \,b_{\cdot,j}(x(t)).
\end{equation}

The following trivial lemma gives a differentiability property of fluid solutions.
\begin{lemma}
\label{lm:differentiable}
Let $x(t)$ be a fluid solution. If $\sum_{i=0}^j x_{\cdot,i}(t)=0$, then $t$ is a point of differentiability if and only if 
\begin{equation}
d \sum_{i=0}^j  (j+1-i) x_{i,\cdot}(t) \le 1.
\end{equation}
\end{lemma}

The following  proposition will be crucial to prove both Theorems~\ref{th1} and~\ref{th2}.
\begin{proposition}
\label{Lyapunov_derivative}
Let $x(t)$ be a fluid solution.
If
$t$ is a point of differentiability, then
\begin{align}
\label{derivative_L_S}
\dot{\mathcal{L}}_S(x(t)) = x_{0,\cdot}(t) - 1 + \lambda.
\end{align}
\end{proposition}
\emph{Proof:}
Let $j^*(t)\bydef \min\{j\ge 0: x_{\cdot,j}(t)>0\} $.
To prove \eqref{derivative_L_S}, we consider the cases $j^*(t)\ge 1$ and $j^*(t) = 0$ separately.
Let us drop the dependency on $t$ for notational simplicity.
%

First, assume that $j^*\ge 1$.
Using Definition~\ref{def1}, we obtain
\begin{subequations}
\label{eq:marginals_horizontal}
\begin{align}
b_{1,\cdot}(x) = & 
x_{2,\cdot}-x_{1,\cdot} +  \lambda d x_{0,\cdot}  -  \mathcal{R}_{0}(x) \frac{x_{1,1}}{x_{\cdot,1}}  \, \mathbf{1}_{\{x_{\cdot,1}>0\}}    - \mathcal{G}_1(x)\\
&+ \sum_{j\ge 2}
   \mathcal{R}_{j-2}(x) \frac{x_{0,j-1}}{x_{\cdot,j-1}} \, \mathbf{1}_{\{ x_{\cdot,j-1} > 0 \}}
 - \mathcal{R}_{j-1}(x) \frac{x_{1,j}}{x_{\cdot,j}}  \,  \mathbf{1}_{\{ x_{\cdot,j} > 0 \}}
 + \mathbf{1}_{\{j=I\}} \mathcal{R}_{I-1}(x) \frac{x_{0,I}}{x_{\cdot,I}} \, \mathbf{1}_{\{x_{\cdot,I}>0\}}\\
b_{i,\cdot}(x) = & x_{i+1,\cdot}   -x_{i,\cdot} 
  -  \mathcal{R}_{i-1}(x) \frac{x_{i,i}}{x_{\cdot,i}}  \, \mathbf{1}_{\{x_{\cdot,i}>0\}} 
  +  \mathcal{R}_{i-2}(x) \frac{x_{i-1,i-1}}{x_{\cdot,i-1}} \, \mathbf{1}_{\{x_{\cdot,i-1}>0\}}  
%
+ \mathcal{G}_{i-1} (x) - \mathcal{G}_{i} (x)\mathbf{1}_{\{i<I\}}\\
+ \sum_{j\ge i+1} &  
   \mathcal{R}_{j-2}(x) \frac{x_{i-1,j-1}}{x_{\cdot,j-1}} \, \mathbf{1}_{\{ x_{\cdot,j-1} > 0 \}}
 - \mathcal{R}_{j-1}(x) \frac{x_{i,j}}{x_{\cdot,j}}  \,  \mathbf{1}_{\{ x_{\cdot,j} > 0 \}}
 + \mathbf{1}_{\{j=I\}} \mathcal{R}_{I-1}(x) \frac{x_{i-1,I}}{x_{\cdot,I}} \, \mathbf{1}_{\{x_{\cdot,I}>0\}},
\end{align}
\end{subequations}
that is,
\begin{subequations}
\label{eq:marginals_horizontal2}
\begin{align}
%
%
%
b_{i,\cdot}(x) = & x_{i+1,\cdot}   -x_{i,\cdot}  -\lambda dx_{i,\cdot}  
 + \mathcal{R}_{j^*-1}(x) \frac{x_{i-1,j^*}}{x_{\cdot,j^*}} 
 - \mathcal{R}_{j^*-1}(x) \frac{x_{i,j^*}}{x_{\cdot,j^*}}, \qquad  \forall i=1,\ldots,j^*-1\\
b_{j^*,\cdot}(x) = & x_{j^*+1,\cdot}   -x_{j^*,\cdot} 
 +  \mathcal{G}_{j^*-1}(x)
 - \mathcal{R}_{j^*-1}(x) \frac{x_{j^*,j^*}}{x_{\cdot,j^*}}
 + \mathcal{R}_{j^*-1}(x) \frac{x_{j^*-1,j^*}}{x_{\cdot,j^*}}\\
b_{j^*+1,\cdot}(x) = & x_{j^*+2,\cdot}   -x_{j^*+1,\cdot}  +  \mathcal{R}_{j^*-1}(x) \frac{x_{j^*,j^*}}{x_{\cdot,j^*}}  \\
b_{i,\cdot}(x) = & x_{i+1,\cdot}   -x_{i,\cdot} , \qquad \forall i\ge j^*+2.
\end{align}
\end{subequations}
Since $t$ is a point of differentiability, then $d \sum_{i=0}^j  (j+1-i) x_{i,\cdot} \le 1$ (by Lemma~\ref{lm:differentiable}) and thus
\begin{equation}
\label{asc878ufcj}
\mathcal{G}_{j^*-1} (x)
=  \lambda d    \sum_{i=0}^{j^*-1} x_{i,\cdot} \quad
\mathcal{R}_{j^*-1} (x) = 
\lambda - \lambda d \sum_{i=0}^{j^*-1}  (j^*-i) x_{i,\cdot}   
\end{equation}
Substituting these expressions in \eqref{eq:Lyapunov_function_derivative}, we get
\begin{subequations}
\label{eq:Lyapunov_1}
\begin{align}
\dot{\mathcal{L}}_S(x)
= 
& \sum_{i\ge 1} i (x_{i+1,\cdot}   -x_{i,\cdot})
+ \sum_{i=1}^{j^*-1}
 i \left( -\lambda dx_{i,\cdot}  
 + \mathcal{R}_{j^*-1}(x) \frac{x_{i-1,j^*}}{x_{\cdot,j^*}} 
 - \mathcal{R}_{j^*-1}(x) \frac{x_{i,j^*}}{x_{\cdot,j^*}}\right) \\
& + 
j^* \left( 
\lambda d    \sum_{i'=0}^{j^*-1} x_{i',\cdot}(t)
 - \mathcal{R}_{j^*-1}(x) \frac{x_{j^*,j^*}}{x_{\cdot,j^*}}
 + \mathcal{R}_{j^*-1}(x) \frac{x_{j^*-1,j^*}}{x_{\cdot,j^*}} \right)\\
&  + (j^*+1) \mathcal{R}_{j^*-1}(x) \frac{x_{j^*,j^*}}{x_{\cdot,j^*}} \\
=  & -\sum_{i\ge 1} x_{i,\cdot}  
+ \lambda d   \sum_{i=0}^{j^*-1} (j^*-i) x_{i,\cdot} 
+ \mathcal{R}_{j^*-1}(x) \sum_{i=1}^{j^*-1}
 i \left(  \frac{x_{i-1,j^*}}{x_{\cdot,j^*}} 
        -  \frac{x_{i,j^*}}{x_{\cdot,j^*}}\right) \\
        &
+ j^*\,\mathcal{R}_{j^*-1}(x) \frac{x_{j^*-1,j^*}}{x_{\cdot,j^*}}
+ \mathcal{R}_{j^*-1}(x) \frac{x_{j^*,j^*}}{x_{\cdot,j^*}} \\
= & 
x_{0,\cdot} - 1 + \lambda d   \sum_{i=0}^{j^*-1} (j^*-i) x_{i,\cdot} + \mathcal{R}_{j^*-1}(x)\\
= & 
x_{0,\cdot} - 1  + \lambda
\end{align}
\end{subequations}
as desired.

Now, assume that $j^* = 0$.
In this case, $x_{0,0}>0$ and the drift~$b(x)$ given in Definition~\ref{def1} takes the linear form
\begin{subequations}
\label{eq:sys_red}
\begin{align} 
\label{eq:sys_red00}
b_{0,0}(x) = &\, -\lambda  + \lambda d (x_{0,\cdot} - x_{0,0} )\\
\label{eq:sys_red01}
b_{0,1}(x) = &\, x_{11} - \lambda d x_{0,1}\\
\label{eq:sys_red11}
b_{1,1}(x) = &\, - x_{11}   + \lambda  + \lambda d (x_{1,\cdot}  - x_{1,1}) \\
b_{0,j}(x) = &\,  x_{1,j} - \lambda d x_{0,j}, \quad j>1 \\
b_{i,j}(x) = &\, x_{i+1,j}-x_{i,j} - \lambda d x_{i,j}, \quad j>i, i\ge 1\\
\label{eq:sys_redii}
b_{i,i}(x) = &\, -  x_{i,i} + \lambda d(x_{i,\cdot} -x_{i,i}) , \quad i>1.
\end{align}
\end{subequations}
By taking summations in \eqref{eq:sys_red}, we obtain
\begin{subequations}
\begin{align} 
b_{1,\cdot}(x) = &\,   \lambda    +   x_{2,\cdot}   -  x_{1,\cdot}    \\
b_{i,\cdot}(x) = &\, x_{i+1,\cdot} \mathbf{1}_{\{i+1\le I\}} - x_{i,\cdot}  \quad \forall i\ge 2.
\end{align}
\end{subequations}
Thus,
\begin{align}
\dot{\mathcal{L}}_S(x) & = \sum_{i\ge 1 } i b_{i,\cdot}(x) = \lambda + \sum_{i\ge 1 } i (x_{i+1,\cdot} \mathbf{1}_{\{i+1\le I\}} - x_{i,\cdot})  = x_{0,\cdot} -1 + \lambda
\end{align}
where in the last equality we have used the normalizing condition.
\hfill$\Box$
\indent

\subsection{Existence and uniqueness of fixed points.}
\label{proof_th1}

Assume that $x$ is a fixed point. Then the system of equations
\begin{equation}
\label{eq:system_without_marginals}
b_{i,j}(x)=0,\quad \forall i,j
\end{equation}
must be satisfied.

Let  $j^*\bydef j^*(x)\bydef \min\{j\ge 0: x_{\cdot,j}>0\} $.
By Definition~\ref{def1}, this means that $\mathcal{R}_j(x)=0$ for all $j\ge j^*$.

If $j^*=0$, 
then $x_{0,0}>0$ and the drift~$b(x)$ given in Definition~\ref{def1} takes the linear form given in \eqref{eq:sys_red}.
Removing one equation from \eqref{eq:system_without_marginals} and adding the normalizing condition $\sum_{i,j}x_{i,j}=1$, we obtain a linear system composed of $(I+1)(I+2)/2$ independent equations and $(I+1)(I+2)/2$ unknowns.
It is easy to see that $x^\star$, where $x_{0,0}^\star=1-\lambda-1/d$, $x_{0,1}^\star=1/d$, $x_{1,1}^\star=\lambda$ and $x_{i,j}^\star=0$ on the remaining coordinates, 
is the unique solution of such system.
It is also clear that $x^\star\in \mathcal{S}$ and $x_{0,0}^\star>0$ if and only if $\lambda < 1-1/d=\lambda_1^*$.
Therefore, in the following we assume that $j^*\ge 1$. 

If \eqref{eq:system_without_marginals} holds true, then we must also have
\begin{equation}
\label{eq:system}
b_{\cdot,j}(x)=0,\quad \forall j.
\end{equation}
Using Definition~\ref{def1}, we obtain
\begin{subequations}
\label{eq:marginals_vertical}
\begin{align}
%
b_{\cdot,1}(x) = &   - \lambda d x_{0,1}   - \mathcal{R}_{0}(x)   \,  \mathbf{1}_{\{ x_{\cdot,1} > 0 \}}   + \lambda d (x_{1,\cdot}-x_{1,1})   +  \lambda -  \mathcal{R}_0(x)
 - \mathcal{G}_1(x)\\
\nonumber
b_{\cdot,j}(x) = & \lambda d x_{j,\cdot}    - \lambda dx_{\cdot,j} - \mathcal{R}_{j-1}(x) \,  \mathbf{1}_{\{ x_{\cdot,j} > 0 \}}  + \mathcal{R}_{j-2}(x) \, \mathbf{1}_{\{ x_{\cdot,j-1} > 0 \}} 
+ \mathcal{G}_{j-1} (x) - \mathcal{G}_{j} (x)  \\
b_{\cdot,I}(x) = & \mathcal{R}_{I-2}(x) \, \mathbf{1}_{\{ x_{\cdot,I-1} > 0 \}}  
+ \mathcal{G}_{I-1} (x)
- \sum_{i=0}^{I-1}\lambda d x_{i,I}.
\end{align}
\end{subequations}
which can be rewritten as
\begin{subequations}
\label{eq:Asuc898dc_C}
\begin{align}
b_{\cdot,j}(x) = & 0, &\forall j =0,\ldots, j^*-1 \\
b_{\cdot,j^*}(x) = & \lambda d x_{j^*,\cdot}    - \lambda dx_{\cdot,j^*} - \mathcal{R}_{j^*-1}(x)  + \mathcal{G}_{j^*-1} (x) \\
\label{eq:Asuc898dc_B}
b_{\cdot,j^*+1}(x) = & \lambda d x_{j^*+1,\cdot}    - \lambda dx_{\cdot,j^*+1} + \mathcal{R}_{j^*-1}(x) \\
b_{\cdot,j}(x) = & \lambda d x_{j,\cdot}    - \lambda dx_{\cdot,j},& \forall j\ge j^*+2.
\end{align}
\end{subequations}
Summing~\eqref{eq:Asuc898dc_C} for all $j\ge j^*+2$, we obtain
\begin{equation}
\nonumber
0 = \sum_{j = j^*+2}^I   x_{j,\cdot} - x_{\cdot,j} 
  = - \sum_{i=0}^{j^*+1} \sum_{j\ge j^*+2} x_{i,j}
\end{equation}
and since $x$ is composed of non-negative components only,
we must have $x_{i,j}=0$ when $i\in\{0,1,\ldots,j^*+1\}$ and $j\in\{j^*+2,\ldots,I\}$.
Using this in \eqref{eq:Asuc898dc_C} when $j=j^*+2$, we obtain that necessarily $x_{j,j'}=0$ for all $j'>j$. Using again this argument when $j=j^*+3$, we obtain  $x_{j,j'}=0$ for all $j'>j^*+3$, and so forth for all $j$.
We have thus shown that the only non-zero elements of $x$ can be on coordinates $(i,j^*)$ and $(i,j^*+1)$, for all $i< j^*+1$, and $(i,i)$ for all $i\ge j^*$. 
Substituting these properties in~\eqref{eq:bij}, we also obtain $b_{j-1,j}(x)=x_{j,j}$, provided that  $j\ge j^*+2$.
On the other hand, \eqref{eq:system_without_marginals} must hold true and therefore $x_{j,j}=0$ for all $j\ge j^*$.
Thus, we have shown that the only non-zero elements of $x$ can be on coordinates 
$(i,j^*)$ and $(i,j^*+1)$, for all $i\le j^*+1$. 
In this case, \eqref{eq:system_without_marginals} simplifies to
\begin{subequations}
\label{scjs8asc}
\begin{align}
0=&b_{i,j^*}(x) =  x_{i+1,j^*}  -  \mathbf{1}_{\{i>0\}} x_{i,j^*} - \lambda d x_{i,j^*}  - \mathcal{R}_{j^*-1}(x) \frac{x_{i,j^*}}{x_{\cdot,j^*}},\qquad\forall i<j^* \\
\label{scjs8ascasvg}
0=&b_{i,j^*+1}(x) = x_{i+1,j^*+1}  -  \mathbf{1}_{\{i>0\}} x_{i,j^*+1} - \lambda d x_{i,j^*+1}  + \mathbf{1}_{\{i>0\}} \mathcal{R}_{j^*-1}(x) \frac{x_{i-1,j^*}}{x_{\cdot,j^*}}, \qquad\forall i<j^*+1 \\
\label{scjsy76utc}
0=&b_{j^*,j^*}(x) =  -x_{j^*,j^*}  -  \mathcal{R}_{j^*-1}(x) \frac{x_{j^*,j^*}}{x_{\cdot,j^*}}  + \lambda d (1-x_{j^*,j^*} -x_{j^*+1,j^*+1}) \\
\label{scjs8}
0=&b_{j^*+1,j^*+1}(x) =  -x_{j^*+1,j^*+1}   +  \mathcal{R}_{j^*-1}(x) \frac{x_{j^*,j^*}}{x_{\cdot,j^*}}
\end{align}
\end{subequations}
where 
$\mathcal{R}_{j^*-1}(x)
= \lambda - \lambda d \sum_{i=0}^{j^*-1}  (j^*-i) (x_{i,j^*}+x_{i,j^*+1})
\ge0$ by Lemma~\ref{lm:differentiable} because a fixed point is a fluid solution that is everywhere differentiable.
In the remainder of the proof, we show that the system in \eqref{scjs8asc} admits a unique solution that satisfies
the normalizing condition $x_{\cdot,j^*}+x_{\cdot,j^*+1}=1$
and $x_{\cdot,j^*}>0$ if and only if $\lambda \in \mathcal{I}_{j^*}$, and that such solution satisfies as well the properties given in the statement of Theorem~\ref{th1}.
This will conclude our proof.

Using that $x_{0,j^*}+x_{0,j^*+1}=1-\lambda$  (by Proposition~\ref{Lyapunov_derivative}), equations \eqref{scjs8asc} imply that
\begin{subequations}
\nonumber
\begin{align}
0=&b_{0,j^*}(x) + b_{1,j^*+1}(x)
= y_{1}   - \lambda d (1-\lambda) - (\lambda d)^2  x_{0,j^*+1}  \\
0=&b_{i,j^*}(x)+b_{i+1,j^*+1}(x) 
 =  x_{i+1,j^*}  +  x_{i+2,j^*+1}  -   (1+ \lambda d) (x_{i,j^*} +x_{i+1,j^*+1}) ,\quad \forall i=1,\ldots,j^*-1 \\
0=&b_{j^*,j^*}(x)+b_{j^*+1,j^*+1}(x) 
%
= \lambda d (1 -x_{j^*,j^*}-x_{j^*+1,j^*+1})  - (x_{j^*,j^*} + x_{j^*+1,j^*+1}).
\end{align}
\end{subequations}
Letting $y_i\bydef x_{i,j^*} + x_{i+1,j^*+1}$, the key observation is that
\begin{subequations}
\begin{align}
\label{eqj76d}
y_1 = &  \lambda d (1-\lambda) + (\lambda d)^2  x_{0,j^*+1}  \\
%
%
%
%
y_{i+1}  = &  (1+ \lambda d) y_{i}, \quad \forall i=1,\ldots,j^*-1 \\
y_{j^*}  = & \lambda d (1 -y_{j^*}),
\end{align}
\end{subequations}
and we notice that the last equation is autonomous.
Previous equations imply that
\begin{subequations}
\begin{align}
%
%
%
%
%
y_{j^*-i}    = &  \frac{y_{j^*}}{(1+ \lambda d)^i},  \quad \forall i=1,\ldots,j^*-1 \\
%
%
\label{sj878f6vy}
y_{j^*}  = & \frac{\lambda d}{1+\lambda d}
\end{align}
\end{subequations}
and using \eqref{eqj76d} we obtain the equation
\begin{equation}
\nonumber
\lambda d (1-\lambda) + (\lambda d)^2  x_{0,j^*+1} = \frac{\lambda d}{(1+ \lambda d)^{j^*}}
\end{equation}
which gives
\begin{equation}
\nonumber
 \lambda d   x_{0,j^*+1} = \lambda-1 + \frac{1}{(1+ \lambda d)^{j^*}}.
\end{equation}
Since it must hold true that $x_{0,j^*}+x_{0,j^*+1}=1-\lambda$ (by Proposition~\ref{Lyapunov_derivative}),  we also obtain
\begin{equation}
\nonumber
\lambda d x_{0,j^*}  = (1+\lambda d  ) (1-\lambda) - \frac{1}{(1+ \lambda d)^{j^*}}.
\end{equation}
Now, in order for such $x$ to be feasible,  we need that both $x_{0,j^*}$ and $x_{0,j^*+1}$ are non-negative.
Using previous expressions, it is not difficult to see that $x_{0,j^*}\ge 0$ and $x_{0,j^*+1}\ge 0$  if and only if
$\lambda \in \mbox{cl}(\mathcal{I}_{j^*})$, where cl$(A)$ denotes the closure of set $A$.

We notice that the normalizing condition can be written as
$1=x_{0,j^*}+x_{0,j^*+1} + x_{1,j^*+1} + \sum_{i=1}^{j^*} y_i $
and using the expressions above it is not difficult to check that it is indeed satisfied.

We now proceed with the construction of the fixed point.
Substituting the properties obtained so far, we notice that
\begin{align}
\nonumber
0  = b_{j^*+1,j^*+1}(x) - b_{j^*,j^*}(x) 
%
%
%
%
 = 2 x_{j^*,j^*}   +  2 \mathcal{R}_{j^*-1}(x) \frac{x_{j^*,j^*}}{x_{\cdot,j^*}}  -  \frac{2\lambda d}{1+\lambda d}
\end{align}
and since necessarily $x_{j^*,j^*}>0$, which is immediately implied by \eqref{scjsy76utc}, we obtain
\begin{align*}
\frac{\mathcal{R}_{j^*-1}(x) }{x_{\cdot,j^*}}  = \frac{1}{x_{j^*,j^*}} \frac{\lambda d}{1+\lambda d} - 1.
\end{align*}
Using this equation in \eqref{scjs8asc} and recalling that $x_{0,j^*}$ has been already explicited,
we get
\begin{subequations}
\begin{align*}
0=&b_{0,j^*}(x) =  x_{1,j^*}  - \lambda d x_{0,j^*}  -  x_{0,j^*} \left(\frac{1}{x_{j^*,j^*}} \frac{\lambda d}{1+\lambda d} - 1\right)\\
%
%
%
%
%
0=&b_{i,j^*}(x) =  x_{i+1,j^*}  -  x_{i,j^*} - \lambda d x_{i,j^*}  -  x_{i,j^*} \left(\frac{1}{x_{j^*,j^*}} \frac{\lambda d}{1+\lambda d} - 1\right) ,~~\forall i=1\ldots,j^*-1
%
%
%
\end{align*}
\end{subequations}
and thus,
\begin{subequations}
\begin{align*}
%
%
x_{1,j^*}  &= \left( \lambda d   -  1 +   \frac{ 1 }{x_{j^*,j^*}} \frac{\lambda d}{1+\lambda d}\right) x_{0,j^*}  \\
%
%
%
%
%
x_{i+1,j^*}  & =  \lambda d \left( 1 +   \frac{1}{x_{j^*,j^*}} \frac{1}{1+\lambda d}\right)  x_{i,j^*} 
=  (\lambda d)^{i} \left( 1 +   \frac{1}{x_{j^*,j^*}} \frac{1}{1+\lambda d}\right)^{i}  x_{1,j^*} 
%
%
%
\end{align*}
\end{subequations}
for all $i=1,\ldots,j^*-1$.
In particular, when $i=j^*-1$, the last equation allows us to identify $x_{j^*,j^*}$ by means of the following polynomial equation
\begin{align}
\nonumber
%
F(x_{j^*,j^*})
\bydef &
(\lambda d)^{j^*-2} \left( 1 +   \frac{1}{x_{j^*,j^*}} \frac{1}{1+\lambda d}\right)^{j^*-1}
\left( \lambda d   -  1 +   \frac{ 1 }{x_{j^*,j^*}} \frac{\lambda d}{1+\lambda d}\right)
\left( (1+\lambda d  ) (1-\lambda) - \frac{1}{(1+ \lambda d)^{j^*}} \right) \\
\label{eq:xjj}
& -x_{j^*,j^*} =0 .
\end{align}
Since $x_{i+1,j^*+1}=y_i-x_{i,j^*}$ and the value of $y_i$ has been already explicited for each $i$,
to conclude the proof of existence and uniqueness of a solution of \eqref{scjs8asc}, it remains to show that previous equation admits a unique root in (0,1] when $\lambda\in \mathcal{I}_{j^*}$.
This property follows easily
once noted that 
$\lim_{x\downarrow 0}F(x)=+\infty$,
$F(1)<0$ if $\lambda \in \mathcal{I}_{j^*}$,
and that $F(x)$ is strictly decreasing if $\lambda \in \mathcal{I}_{j^*}$.



\subsection{Bounds on fluid mass.}
\label{section:LS_LM}

Whenever $x(t)$ is differentiable, \eqref{asc878ufcj} and~\eqref{eq:Asuc898dc_C} imply
\begin{subequations}
\label{eq:Lyapunov_1}
\begin{align}
\dot{\mathcal{L}}_M
=&  \sum_{j=1}^I j \,b_{\cdot,j}(x) \\
=&  j^* \left( \lambda d x_{j^*,\cdot}    - \lambda dx_{\cdot,j^*} - \mathcal{R}_{j^*-1}(x)  + \mathcal{G}_{j^*-1} (x) \right) \\
& + (j^*+1) \left( \lambda d x_{j^*+1,\cdot}    - \lambda dx_{\cdot,j^*+1} + \mathcal{R}_{j^*-1}(x) \right) \\
& + \lambda d \sum_{j\ge j*+2} j  ( x_{j,\cdot}    - x_{\cdot,j} ) \\
=& \lambda + \lambda d \sum_{i=0}^{j^*-1}  i x_{i,\cdot} +   \lambda d \sum_{j\ge j^*} j  ( x_{j,\cdot}    - x_{\cdot,j} ) \\
=& \lambda + \lambda d \, \mathcal{L}_S - \lambda d \, \mathcal{L}_M.
\end{align}
\end{subequations}


In a fixed point $x$, we must have $\dot{\mathcal{L}}_M=0$. This condition gives \eqref{skcj878dc} and
since $x_{\cdot,j^\star}+x_{\cdot,j^\star+1}=1$ (by Theorem~\ref{th1}) we obtain
\begin{align*}
\mathcal{L}_S(x) = \mathcal{L}_M(x) -  \frac{1}{d} &\,  = j^\star x_{\cdot,j^\star} + (j^\star+1) x_{\cdot,j^\star+1} -  \frac{1}{d} 
%
%
 = j^\star  +  x_{\cdot,j^\star+1} -  \frac{1}{d}
\end{align*}
and \eqref{skcj878dcsvdf} holds true because $x_{\cdot,j^\star+1}^\star\in [0,1]$.

\subsection{Global stability.}
\label{proof_th2}

To prove Theorem~\ref{th2}, we first introduce the following lemma.

\begin{lemma}
\label{lemma_x01}
Assume $\lambda < 1-\frac{1}{d}$.
Let $x(t)$ be a fluid solution such that  $x_{0,0}(t)>0$ for all $t\ge 0$. Then,  \eqref{eqsjauis} holds true.
\end{lemma}
\emph{Proof:}
If $x_{0,0}>0$, then the drift $b(x)$ has the linear form given in~\eqref{eq:sys_red}.
Thus, if $x_{0,0}(t)>0$ for all $t\ge 0$, then the fluid solution $x(t)$ is uniquely determined by the ODE system 
\begin{subequations}
\label{eq:sys_red_asc}
\begin{align} 
\label{eq:sys_red00_asc}
\dot{x}_{0,0} = &\, -\lambda  + \lambda d (x_{0,\cdot} - x_{0,0} )\\
\label{eq:sys_red01_asc}
\dot{x}_{0,1} = &\, x_{11} - \lambda d x_{0,1}\\
\dot{x}_{1,1} = &\, - x_{11}   + \lambda  + \lambda d (x_{1,\cdot}  - x_{1,1}) \\
\dot{x}_{0,j} = &\,  x_{1,j} - \lambda d x_{0,j}, \quad j>1 \\
\dot{x}_{i,j} = &\, x_{i+1,j}-x_{i,j} - \lambda d x_{i,j}, \quad j>i, i\ge 1\\
\label{eq:sys_redii_asc}
\dot{x}_{i,i} = &\, -  x_{i,i} + \lambda d(x_{i,\cdot} -x_{i,i}) , \quad i>1.
\end{align}
\end{subequations}
We have already shown in Section~\ref{proof_th1} that $x^\star$ is the unique fixed point of such linear system.
The equations \eqref{eq:sys_red01_asc}-\eqref{eq:sys_redii_asc} do not depend on $x_{0,0}$ and thus form an autonomous ODE system.
This means that we can safely remove the equation~\eqref{eq:sys_red00_asc} and recall that~$x_{0,0}(t)$ can be uniquely obtained by using the normalizing condition, i.e., $ x_{0,0}(t)= 1 - \sum_{(i,j)\neq(0,0)} x_{i,j}(t)$.
%
The ODE system \eqref{eq:sys_red01_asc}-\eqref{eq:sys_redii_asc} has the linear form $\dot{x}=A x + p$ where
$A$ is a triangular matrix and $p$ is a column vector, 
and it is clear that the eigenvalues of $A$ are $-1, -\lambda d $ and $-(1+\lambda d)$.
%
%
Since the eigenvalues of $A$ are strictly negative,
it follows from standard results in ODE theory that
$x(t) = x^\star + e^{At} (x(0)-x^\star)$.
Thus, \eqref{eqsjauis} follows by the norm bound on the exponential matrix. 
\hfill$\Box$
\indent

For the fluid solution $x(t)$, either $x_{0,0}(t)>0$  for all $t\ge 0$, in which case Theorem~\ref{th2} follows directly by previous lemma, or 
$x_{0,0}(t_0)=0$ for some~$t_0$, that is the case we study in the following.
Without loss of generality, let us assume $t_0>0$.

If $t_0$ is a point of differentiability of $x_{0,0}(\cdot)$, then Lemma~\ref{lm:differentiable} implies that $d x_{0,\cdot}(t_0) \le 1$.
If $t_0$ is not a point of differentiability of $x_{0,0}(\cdot)$, then we still have $d x_{0,\cdot}(t_0) \le 1$ because
either
there exists $\delta$ such that 
$x_{0,0}(t)=0$ for all  $t\in[t_0-\delta,t_0]$ and the inequality holds true again by Lemma~\ref{lm:differentiable}
or 
there exists a sequence $t_n\uparrow t_0$, $n\ge 1$, such that $t_n<t_{n+1}<t_0$ where $x_{0,0}(t_n)>0$ and $0>\dot{x}_{0,0}(t_n)= d \lambda  (x_{0,\cdot}(t_n) - x_{0,0}(t_n)) - \lambda$ (by \eqref{eq:b00}) for all $n$, which implies $d \lambda  x_{0,\cdot}(t_0)  - \lambda = \lim_{n\to\infty}d \lambda  (x_{0,\cdot}(t_n) - x_{0,0}(t_n))  - \lambda \le 0$.

Lemma~\ref{lm:differentiable} also ensures that $d x_{0,\cdot}(t) \le 1$ on $[t_0,\infty)$ as long as $x_{0,0}(t)=0$.
Substituting $d x_{0,\cdot}(t) \le 1$ in \eqref{derivative_L_S}, we obtain
$\dot{\mathcal{L}}_S(x(t))  \le   1/d - 1 + \lambda$.
Since $\lambda<1-1/d$ by hypothesis, this means that $x_{0,0}(t)$ cannot remain equal to zero on $[t_0,\infty)$ because 
$\mathcal{L}_S(x(t))$ would be decreasing in $t$ with derivative bounded away from zero and necessarily
$\mathcal{L}_S(x)\ge 0$ for all~$x$.
This implies that
$t^*\bydef \inf\{t\ge t_0 : x_{0,0}(t)>0\}<\infty$ must exist.
Since $x(t)$ is continuous and $b(x)$ is linear when $x_{0,0}>0$, there exists $\delta>0$ such that  $x_{0,0}(t)$ is both positive and increasing on $(t^*,t^*+\delta]$.
This means that $0< \dot{x}_{0,0}(t) = \lambda d(x_{0,\cdot}(t)- x_{0,0}(t)) - \lambda $ for all $t\in (t^*,t^*+\delta]$ (by \eqref{eq:b00}), which implies 
$\lim_{t\downarrow t^*} x_{0,\cdot}(t)- x_{0,0}(t) =x_{0,\cdot}(t^*)\ge  1/d$.
On the other hand, on a left neighborhood of~$t^*$, 
$x_{0,\cdot}(t^*) = \lim_{t\uparrow t^*} x_{0,\cdot}(t)\le  1/d$,
and thus
$x_{0,\cdot}(t^*)=1/d$.
Summarizing,  we have obtained
\begin{subequations}
\label{eq:uas8c1}
\begin{align}
& x_{0,\cdot}(t^*)=  \tfrac{1}{d},
\quad
x_{0,0}(t^*)=0 & \\
&
x_{0,0}(t)>0,
\quad
x_{0,\cdot}(t)-x_{0,0}(t)>1/d, &
\forall t\in(t^*,t^*+\delta].
\end{align}
\end{subequations}
These conditions, together with the fact that $b(x)$ is linear when $x_{0,0}>0$, imply that
also the function 
$w_0(t) \bydef x_{0,\cdot}(t)-x_{0,0}(t)$ must be increasing on a right neighborhood of $t^*$. 
Thus,
\begin{equation}
\label{eq:uas8c2}
\dot{w}_{0}(t^*) = \lim_{t\downarrow t^*} \dot{w}_{0}(t)\ge 0.
\end{equation}
On $[t^*,t^*+\delta]$, we have shown that the fluid solution $x(t)$ is uniquely determined by the solution of the ODE system \eqref{eq:sys_red_asc} where the initial condition $x(t^*)$ is such that 
\eqref{eq:uas8c1} and \eqref{eq:uas8c2} hold true.
In the remaining part of the proof, we study \eqref{eq:sys_red_asc} under these conditions and show that $x_{0,0}(t)>0$ for all $t>t^*$. This will conclude the proof in view of Lemma~\ref{lemma_x01}.

Without loss of generality and by means of a time shift, let us assume $t^*=0$.
By taking proper summations in \eqref{eq:sys_red_asc}, we obtain
\begin{align}
\label{w0_dvfc}
\dot{w}_{0,\cdot} = &\,  -\lambda d w_0 + x_{1,\cdot} 
\end{align}
Now, given that $w_0(0)=1/d$, \eqref{eq:uas8c2} ensures that 
$\dot{w}_{0}(0) = -\lambda d w_0(0) + x_{1,\cdot}(0) \ge 0 $, which means
\begin{equation}
\label{eq:x1_gt_lambda}
x_{1,\cdot}(0)\ge \lambda.
\end{equation}
By taking proper summations in \eqref{eq:sys_red_asc}, we also obtain
\begin{subequations}
\label{wfiu8e7d_sys}
\begin{align} 
\label{wfiu8e7d_1}
\dot{x}_{0,\cdot} = &\,  -\lambda + x_{1,\cdot}    \\
\label{wfiu8e7d_2}
\dot{x}_{1,\cdot} = &\,   \lambda    +   x_{2,\cdot}   -  x_{1,\cdot}    \\
\dot{x}_{i,\cdot} = &\, x_{i+1,\cdot} \mathbf{1}_{\{i+1\le I\}} - x_{i,\cdot}  \quad \forall i\ge 2.
\end{align}
\end{subequations}
and solving for such autonomous ODE system,
\begin{align*} 
x_{I,\cdot}(t) = & \, x_{I,I} =   x_{I,I}(0) \, e^{-t} \\
x_{I-1,\cdot}(t) = & \, x_{I,I}(0) \, e^{-t} - x_{I-1,\cdot} =  (x_{I,I}(0) t + x_{I-1,\cdot}(0) )\,e^{-t}  \\
%
%
%
x_{I-i,\cdot}(t) = & \,  e^{-t}  \sum_{j=0}^i x_{I-i+j,\cdot}(0) \, \frac{t^j}{j!} 
\end{align*}
for all $i=2,\ldots, I-2$, and thus
\begin{equation*}
x_{2,\cdot}(t) =   e^{-t}  \sum_{j=0}^{I-2} x_{2+j,\cdot}(0) \, \frac{t^j}{j!} .
\end{equation*}
Substituting previous equation in \eqref{wfiu8e7d_2}, we obtain
\begin{equation}
\label{x1_ascds}
\dot{x}_{1,\cdot} =   \lambda -  x_{1,\cdot}  + e^{-t}  \sum_{j=0}^{I-2} x_{2+j,\cdot}(0) \, \frac{t^j}{j!} 
\end{equation}
and solving for such ODE we obtain,
\begin{equation*}
x_{1,\cdot}(t) =   \lambda   + (x_{1,\cdot}(0)-\lambda) e^{-t} +  e^{-t}  \sum_{j=1}^{I-1} x_{1+j,\cdot}(0) \, \frac{t^j}{j!} 
\end{equation*}
Substituting previous equation in \eqref{wfiu8e7d_1}, we obtain
\begin{align*}
\dot{x}_{0,\cdot}(t) = &\,  (x_{1,\cdot}(0)-\lambda) e^{-t} +  e^{-t}  \sum_{j=1}^{I-1} x_{1+j,\cdot}(0) \, \frac{t^j}{j!}.
\end{align*}
Integrating both sides, we obtain
\begin{align*}
x_{0,\cdot}(t)-\frac{1}{d}
= &\, x_{1,\cdot}(0)-\lambda - (x_{1,\cdot}(0)-\lambda) e^{-t}  +   \sum_{j=1}^{I-1} \frac{x_{1+j,\cdot}(0)}{j!} \, \int_0^t e^{-s} s^j  \mbox{d}s.
\end{align*}
Similarly, substituting \eqref{x1_ascds} in \eqref{w0_dvfc} and solving for $w_0(t)$, we obtain
\begin{align*}
w_0(t) 
%
= &\,   
\frac{1}{d}   +   \frac{x_{1,\cdot}(0) - \lambda}{\lambda d -1} \left( e^{-t} - e^{-\lambda d t} \right) 
+ e^{-\lambda d t}  \sum_{j=1}^{I-1}  x_{1+j,\cdot}(0)   \int_{0}^t \frac{s^j}{j!}  e^{(\lambda d -1)s} \mbox{d}s
\end{align*}
and thus
\begin{align*}
& x_{0,0}(t) =  \, x_{0,\cdot}(t) -w_{0}(t)  \\
%
&\,   = 
 (x_{1,\cdot}(0)-\lambda) (1-e^{-t} )
-   \frac{x_{1,\cdot}(0) - \lambda}{\lambda d -1} \left( e^{-t} - e^{-\lambda d t} \right)  
+   \sum_{j=1}^{I-1} x_{1+j,\cdot}(0)  \, \int_0^t \frac{s^j}{j!}
e^{ -s} \left(  1   -   e^{ -\lambda d (t-s)} \right)
\mbox{d}s.
\end{align*}
We now use the condition~\eqref{eq:x1_gt_lambda}.
If $x_{1,\cdot}(0) = \lambda$, then $\sum_{i\ge 2} x_{i,\cdot}(0) = 1-1/d-\lambda >0 $ and therefore
\begin{align*}
x_{0,0}(t) =& 
\sum_{j=1}^{I-1} \frac{x_{1+j,\cdot}(0)}{j!} \, \int_0^t s^j
e^{ -s} \left(  1   -   e^{ -\lambda d (t-s)} \right)
\mbox{d}s
\ge 
\left( 1-\frac{1}{d}-\lambda\right)  \min_{j=1}^{I-1}  \int_0^t \frac{s^j}{j!}
e^{ -s} \left(  1   -   e^{ -\lambda d (t-s)} \right)
\mbox{d}s >0
\end{align*}
for all $t>0$, as desired.
If $x_{1,\cdot}(0) > \lambda$, then 
\begin{align*}
x_{0,0}(t) \ge  &\, 
 (x_{1,\cdot}(0)-\lambda) (1-e^{-t} )
-   \frac{x_{1,\cdot}(0) - \lambda}{\lambda d -1} \left( e^{-t} - e^{-\lambda d t} \right) .
\end{align*}
Given that $x_{0,0}(0)=0$, 
to conclude that $x_{0,0}(t)>0$ for all $t>t^*=0$ it is sufficient to show that the RHS of last equation is strictly increasing in $t$.
This follows easily once noted that the derivative of the RHS of last equation is strictly positive if and only if
\begin{align*}
\frac{e^{-t}   -  e^{-\lambda d t}}{\lambda d -1}   > 0,\quad \forall t>0.
\end{align*}

\section{Conclusions}
\label{sec:con}

In this paper, we have provided new insights on randomized load balancing:
if a load balancer is endowed with a local memory storing the last observation collected on each server, the celebrated power-of-$d$-choices algorithm can be made asymptotically optimal in the sense that arriving jobs can be always routed to idle servers.
Our approach provides an algorithm that is \emph{both} fluid ($N\to\infty$) and heavy-traffic ($\lambda\uparrow 1$) optimal while employing a fair control message rate that scales linearly with the system size~$N$.
This means that randomized load balancing can be made robust to orthogonal variations of both $N$ and $\lambda$, which can for instance occur in presence of unexpected workload peaks or server breakdowns.

On the practical side, Algorithm~\ref{alg} can be improved in several ways to enhance performance: 
\begin{itemize}
 \item Server selections can be made without replacement, instead of with replacement. In view of the results in Gast and Van Houdt \cite{Gast2017}, this may also improve the convergence speed of $X^N$ to fluid solutions.

 \item Since the action of sampling $(0,0)$-servers does not bring any additional information to the load balancer,
 server selections can be restricted to $(\cdot,j)$-servers, with $j\ge 1$.
 
 \item Upon a job arrival, if $i$ is both the least load of the $d$ sampled servers and the least observation contained in the memory immediately before the last sampling, then the job may be randomly assigned to one of the $(\cdot,i)$-server known to the load balancer immediately before the sampling. In fact, by the time of the last update, one of such servers may have decreased its load.
 
 \item At any point in time, it is clear that the observation collected on a specific server is an upper bound on the actual state of that server.
This observation leads us to consider a variant of SQ$(d,N)$ where to each server is associated a \emph{timer} representing the age of its observation.
Specifically, at the moment where a new job arrives, all timers are incremented by one except the ones associated to the $d$ sampled servers, whose timers are set to zero.
Then, the load-balancer dispatches the job to a server with the lowest stored observation and for which the timer is the largest.

\end{itemize}

It is intuitive that all the above variations of SQ$(d,N)$ yield performance improvements.
It may be less intuitive that at the fluid scale only the last variation can yield performance improvements.
In view of the results presented in this paper, such improvements can only appear when $\lambda>1-1/d$.
We leave this subject as future research. 

A last variant of Algorithm~\ref{alg} consists in swapping Lines 4--7 and 8--10.
This is meant to perform each job assignment before the $d$ servers are sampled (the collected information will be thus used for future assignments).
It can be easily shown, mutatis mutandis, that this yields the same fluid limit.
This is not surprising: if $x_{0,0}>0$, the number of zeros in the memory is proportional to $N$ and getting $d$ more observations does not make $x_{0,0}$ zero.


In our analysis, we have assumed that each server has a finite buffer of size~$I$. We conjecture that our results generalize to the case where $I=\infty$ and a first step to prove this claim consists in adapting the proofs of Lemma~\ref{LM:tightness} and Proposition~\ref{LM:tightness} on coordinates $(i,i)$ only.
Provided that servers are initially empty, this conjecture is coherent with the numerical observation that the Lyapunov function~$\mathcal{L}_S(x(t))$ monotonically increases in $t$ to its limit point, which is necessarily less than~$j^\star+1$; see Figure~\ref{fig:simulations}.
If $\lambda<1-1/d$ and $x_{0,0}(0)>0$, this can be easily proven by using Lemma~\ref{lemma_x01}, which ensures that the drift function~$b$ takes the linear form in~\eqref{eq:sys_red}. When $\lambda\ge 1-1/d$, a proof is complicated by the involved structure of $\mathcal{L}_S$.


Our model can be generalized to a setting where servers have bin-packing constraints.
Specifically, each server has $B$ units of a resource, there are $R$ types of jobs, type-$r$ jobs requires $b_r$ units of resource, and an arriving job is `blocked' if it does not find the required amount of resource at the server.
Memoryless power-of-$d$-choice strategies have been recently applied to this type of models in Xie et al. \cite{Xie2015_PdC_binpacking}, though the resulting blocking probability does not converge to zero in the fluid limit. 
A further direction for future research aims at evaluating whether or not a local memory at the dispatcher can still be exploited to achieve fluid optimality in this setting.


%
%
%
\appendix

\section{Proof of Lemma~\ref{lemma_ascasji1}}

We give a proof when $\epsilon\downarrow 0$ as the same arguments can be applied when $\epsilon\uparrow 0$.
Let $s_{i,j}\bydef s_{i,j}(\overline{x}) \bydef  \sum_{i'=0}^{i-1}  \overline{x}_{i',\cdot} + \sum_{j'\ge i}^j  \overline{x}_{i,j'}$.
For any $\epsilon>0$ and all $k$ sufficiently large, Lemma~\ref{lemma:X_close_x} states that the inclusions
\begin{subequations}
\begin{align*}
( s_{i,i}(t)+2 \epsilon L (I+1)^2,  s_{i,\cdot}(t) -2 \epsilon L (I+1)^2   ] 
& \subseteq ( S_{i,i}^{N_k}(t_n^{{N_k},\lambda-}), S_{i,\cdot}^{N_k}(t_n^{{N_k},\lambda-})   ] \\
&  \subseteq  ( s_{i,i}(t)-2 \epsilon L (I+1)^2, s_{i,\cdot}(t) +2 \epsilon L (I+1)^2  ]
\end{align*}
\end{subequations}
hold true, 
as we recall that under the coupled construction given in Section~\ref{sec:constructionX} we have that $t_n^{N_k,\lambda-}\in (t,t+\epsilon]$ when $n \in \{ \mathcal{N}_\lambda({N_k} t)+1, \ldots, \mathcal{N}_\lambda({N_k} (t+\epsilon))\}$.
Using these inclusions and Lemma~\ref{GC_F}, we obtain
\begin{subequations}
\begin{align*} 
 - \lambda \epsilon + d \lambda  \left(\overline{x}_{i,\cdot}(t) - \overline{x}_{i,i}(t) - 4 \epsilon L (I+1)^2 \right) \epsilon  
 &\le 
\lim_{k\to\infty}  \frac{1}{N_k} \sum_{n=\mathcal{N}_\lambda(N_k t)+1}^{\mathcal{N}_\lambda(N_k (t+\epsilon))}
\sum_{p=1}^d  \mathbb{I}_{ ( S_{i,i}^{N_k}(t_n^{N_k,\lambda-}),  S_{i,I}^{N_k}(t_n^{N_k,\lambda-})   ] }^{(V_n^p)}  \\
%
& \le   - \lambda \epsilon + d \lambda  \left(\overline{x}_{i,\cdot}(t) - \overline{x}_{i,i}(t) + 4 \epsilon L (I+1)^2 \right) \epsilon
\end{align*}
\end{subequations}
for all $\epsilon>0$, and
dividing these inequalities by $\epsilon$ and letting $\epsilon\downarrow 0$, we obtain~\eqref{lemma_ascasji12AAA}.

Finally, \eqref{lemcj_hcyty} is proven using the same argument.

\section{Proof of Lemma~\ref{lemma:R}}

At the beginning of the proof of Lemma~\ref{prop_R0}, we have already shown that 
$R_j(t)=0$ if  $\sum_{j'=0}^j \overline{x}_{\cdot,j'}(t)>0$.
Thus, in the following we assume that $\sum_{j'=0}^j \overline{x}_{\cdot,j'}(t)=0$.

First, for the indicator function in \eqref{eq:scuas78cas3}, we will use that
\begin{subequations}
\label{ascuja7sc8as2}
\begin{align}
\mathbf{1}_{\{X_{0,0}^N(t_n^{N,\lambda-}) + \sum_{j'=1}^I M_{0,j',n}^N>0\}}
= &  1 - \mathbf{1}_{\{X_{0,0}^N(t_n^{N,\lambda-}) = 0\}} \mathbf{1}_{\{\sum_{j'=1}^I M_{0,j',n}^N=0\}} \\
= &  1 - \mathbf{1}_{\{X_{0,0}^N(t_n^{N,\lambda-}) = 0\}} \prod\limits_{p=1}^d \mathbb{I}_{ ( 1-Z_{1}^N(t_n^{N,\lambda-}), 1 ] }^{(V_n^p)}.
\end{align}
\end{subequations}
We also notice the following sequence of equalities (see \eqref{def:F_ijn} and \eqref{def:F_iin} for the definition of $F_{i,j,n}^N$)
\begin{subequations}
\label{ascuja7sc8as}
\begin{align}
& \sum_{i=0}^j F_{i,j,n}^N 
= 
\mathbb{I}_{
\left( 
0,
X_{\cdot,j}^N(t_n^{N,\lambda-}) + \overline{M}_{j,n}^N - \underline{M}_{j,n}^N
\right]
}^{
(W_n (X_{\cdot,j}^N(t_n^{N,\lambda-}) + \overline{M}_{j,n}^N - \underline{M}_{j,n}^N )  )}\\
 &=   \mathbf{1}_{   \{ X_{\cdot,j}^N(t_n^{N,\lambda-}) + \overline{M}_{j,n}^N - \underline{M}_{j,n}^N >0 \}  } 
 =    1 -  \mathbf{1}_{   \{ X_{\cdot,j}^N(t_n^{N,\lambda-}) - \underline{M}_{j,n}^N  =0 \}}  \mathbf{1}_{   \{ \overline{M}_{j,n}^N =0 \}  } \\
 &=   1 -  \mathbf{1}_{   \{ X_{\cdot,j}^N(t_n^{N,\lambda-}) - \underline{M}_{j,n}^N  =0 \}}  \prod_{p=1}^d \mathbb{I}_{ ( 0, 1-Z_{j}^{N}(t_n^{N,\lambda -}) ] \bigcup\, (1 - Z_{j+1}^{N}(t_n^{N,\lambda -}) ,1]} ^{(V_n^p)}\\
\nonumber
&=   1 -
\left( 
\mathbf{1}_{   \{ X_{\cdot,j}^N(t_n^{N,\lambda-}) =0 \}}
+
\sum_{p=1}^d  \mathbf{1}_{   \{ NX_{\cdot,j}^N(t_n^{N,\lambda-})  =p \}} \mathbf{1}_{   \{ N\underline{M}_{j,n}^N  =p \}}
\right)
\prod_{p=1}^d \mathbb{I}_{ ( 0, 1-Z_{j}^{N}(t_n^{N,\lambda -}) ] \bigcup\, (1 - Z_{j+1}^{N}(t_n^{N,\lambda -}) ,1]}^{(V_n^p)}\\
\label{fjw8tg8irfasc}
 &=   1 -
\mathbf{1}_{   \{ X_{\cdot,j}^N(t_n^{N,\lambda-}) =0 \}}
\prod_{p=1}^d \mathbb{I}_{ ( 0, 1-Z_{j}^{N}(t_n^{N,\lambda -}) ] \bigcup\, (1 - Z_{j+1}^{N}(t_n^{N,\lambda -}) ,1]}^{(V_n^p)}\\
 \quad &
\label{fjw8tg8irf}
 -\sum_{p=1}^d  \mathbf{1}_{   \{ NX_{\cdot,j}^N(t_n^{N,\lambda-})  =p \}} \mathbf{1}_{   \{ N\underline{M}_{j,n}^N  =p \}}
\times
\prod_{p=1}^d \mathbb{I}_{ ( 0, 1-Z_{j}^{N}(t_n^{N,\lambda -}) ] \bigcup\, (1 - Z_{j+1}^{N}(t_n^{N,\lambda -}) ,1]}^{(V_n^p)}.
%
%
%
%
\end{align}
\end{subequations}
For the summation terms in \eqref{fjw8tg8irf}, we observe that
\begin{align}
\nonumber
&\mathbf{1}_{   \{ N\underline{M}_{j,n}^N  =p \}}
=  \mathbf{1}_{   \left\{ \sum_{i=0}^j N M_{i,j,n}^N  =p \right \}} 
= \mathbf{1}_{   \left \{ \sum\limits_{i=0}^j \sum\limits_{q=1}^d  \mathbb{I}_{ ( S_{i,j-1}^N(t_n^{N,\lambda-}),  S_{i,j}^N(t_n^{N,\lambda-})   ] }^{(V_n^q)}  =p \right \}} \\
\nonumber
 &= \mathbf{1}_{   \left \{ \sum\limits_{q=1}^d  \mathbb{I}_{ \bigcup_{i=0}^j ( S_{i,j-1}^N(t_n^{N,\lambda-}),  S_{i,j}^N(t_n^{N,\lambda-})   ] }^{(V_n^q)}  =p \right \}}\\
%
\nonumber
 &=
\sum_{I\subseteq \{1,\ldots,d\}: \|I\|=p} 
\prod_{q\in I}  \mathbb{I}_{ \bigcup_{i=0}^j ( S_{i,j-1}^N(t_n^{N,\lambda-}),  S_{i,j}^N(t_n^{N,\lambda-})   ] }^{(V_n^q)}
\times
\prod_{q\notin I}  \mathbb{I}_{ \neg \bigcup_{i=0}^j ( S_{i,j-1}^N(t_n^{N,\lambda-}),  S_{i,j}^N(t_n^{N,\lambda-})   ] }^{(V_n^q)}\\
\nonumber
 &\le
\sum_{I\subseteq \{1,\ldots,d\}: \|I\|=p} 
\prod_{q\in I}  \mathbb{I}_{ \bigcup_{i=0}^j ( S_{i,j-1}^N(t_n^{N,\lambda-}),  S_{i,j}^N(t_n^{N,\lambda-})   ] }^{(V_n^q)}
\end{align}
where $\neg A$ denotes the complement of set $A$ and we have defined
$S_{j,j-1}^N(t_n^{N,\lambda-}) \bydef S_{j-1,I}^N(t_n^{N,\lambda-})$.
Thus, for any $\epsilon>0$, on the interval $[t,t+\epsilon]$,
\eqref{X_close_x} ensures that
\begin{align}
\label{askciudfg}
%
0 \le
\mathbf{1}_{   \{ N_k \underline{M}_{j,n}^{N_k }  =p \}} \le 
%
\sum_{
{\scriptstyle I\subseteq \{1,\ldots,d\}:} \atop {\scriptstyle \|I\|=p}  
}
\prod_{q\in I}  \mathbb{I}_{ \bigcup_{i=0}^j ( 0\vee s_{i,j-1}- \epsilon,  s_{i,j} + \epsilon   ] }^{(V_n^q)}
\end{align}
for all $k$ sufficiently large,
where $s_{i,j}\bydef s_{i,j}(\overline{x}(t))\bydef \sum_{i'=0}^{i-1}  \overline{x}_{i',\cdot}(t) + \sum_{j'\ge i}^j  \overline{x}_{i,j'}(t)$ for all $i\le j$ and $s_{j,j-1} = s_{j-1,I}$.

%
Let us treat the cases $j=1$ and $j>1$ separately.

Assume for now $j=1$. Substituting \eqref{ascuja7sc8as2} in the sample path expressions~\eqref{eq:shucyu} and~\eqref{eq:shucyuascs}, we obtain
\begin{subequations}
\nonumber
\begin{align}
X_{\cdot,1}^N(t)  =  & X_{\cdot,1}^N(0 )  + \frac{1}{N} \sum_{n=1}^{\mathcal{N}_\lambda(N t)}
\sum_{p=1}^d  \mathbb{I}_{ ( S_{1,1}^N(t_n^{N,\lambda-})  , S_{1,I}^N(t_n^{N,\lambda-})  ] }^{(V_n^p)}
\\
& + \frac{1}{N} \sum_{n=1}^{\mathcal{N}_\lambda(N t)}
\left( 1 - \mathbf{1}_{\{R_{0}^N(t_n^{N,\lambda-}) = 0\}} \prod\limits_{p=1}^d \mathbb{I}_{ ( 1-Z_{1}^N(t_n^{N,\lambda-}), 1 ] }^{(V_n^p)} \right)\\
&  - \frac{1}{N} \sum_{n=1}^{\mathcal{N}_\lambda(N t)}
\sum_{p=1}^d  \mathbb{I}_{ ( S_{0,1}^N(t_n^{N,\lambda-})-X_{0,1}^N(t_n^{N,1-}),  S_{0,1}^N(t_n^{N,\lambda-})   ] }^{(V_n^p)} \\
%
%
\label{eq:ajicas94}
&  - \frac{1}{N} \sum_{n=1}^{\mathcal{N}_\lambda(N t)}  \mathbf{1}_{ \left\{R_{0}^N(t_n^{N,\lambda-}) = 0 \right\} } 
\left (  F_{0,1,n}^{N} + F_{1,1,n}^{N}  \right) 
\prod\limits_{p=1}^d \mathbb{I}_{ ( 1-Z_1^N(t_n^{N,\lambda-}), 1 ] }^{(V_n^p)}.
%
%
%
%
\end{align}
\end{subequations}
Since $\overline{x}_{0,0}(t) + \overline{x}_{\cdot,1}(t)=0$ and $t$ is a point of differentiability, necessarily $\dot{\overline{x}}_{\cdot,1}(t)=0$ and thus
\begin{subequations}
 \label{ascjh87tfasc}
\begin{align}
0 = & \dot{\overline{x}}_{\cdot,1}(t) \\
 = & \lim_{\epsilon\to 0} \frac{1}{\epsilon}\lim_{k\to\infty}  X_{\cdot,1}^{N_k}(t+\epsilon)-X_{\cdot,1}^{N_k}(t) \\
 \label{ascjh87tf_asc}
 = & \lambda d \overline{x}_{1,\cdot}(t)
 + \lambda -  R_0(t)
\\
& 
 \label{ascjh87tf}
 -
\lim_{\epsilon\to 0} \frac{1}{\epsilon}\lim_{k\to\infty}
 \frac{1}{N_k} 
 \sum_{n=\mathcal{N}_\lambda(N_k t)+1}^{\mathcal{N}_\lambda(N_k (t+\epsilon))} 
%
\mathbf{1}_{ \left\{R_{0}^{N_k}(t_n^{N_k,\lambda-}) = 0 \right\} } 
\left (  F_{0,1,n}^{N_k} + F_{1,1,n}^{N_k}  \right) 
\prod\limits_{p=1}^d \mathbb{I}_{ ( 1-Z_1^{N_k}(t_n^{N_k,\lambda-}), 1 ] }^{(V_n^p)}.
\end{align}
\end{subequations}
The terms in \eqref{ascjh87tf_asc} are a direct application of Lemmas~\ref{GC_F} and \ref{lemma_ascasji1}.
Equation \eqref{ascjh87tfasc} also states that the limit in \eqref{ascjh87tf} exists.
Now, using \eqref{ascuja7sc8as}-\eqref{askciudfg} and \eqref{X_close_x}, for any $\epsilon>0$ 
\begin{subequations}
\label{jvr7h6tg}
\begin{align}
& 1 - \mathbf{1}_{   \{ R_{j}^{N_k}(t_n^{N_k,\lambda-}) =0 \}}
\prod_{p=1}^d \mathbb{I}_{ ( 0, 1-Z_{j}^{N_k}(t_n^{N_k,\lambda -}) ] \bigcup\, (1 - Z_{j+1}^{N_k}(t_n^{N_k,\lambda -}) ,1]}^{(V_n^p)}  \\
& \ge \sum_{i=0}^j  F_{i,j,n}^{N_k}   \ge 
\\
 &
 \label{sjc8asc9h}
%
%
%
%
%
%
1 - \mathbf{1}_{   \{ R_{j}^{N_k}(t_n^{N_k,\lambda-}) =0 \}}
\prod_{p=1}^d \mathbb{I}_{ ( 0, 1-Z_{j}^{N_k}(t_n^{N_k,\lambda -}) ] \bigcup\, (1 - Z_{j+1}^{N_k}(t_n^{N_k,\lambda -}) ,1]}^{(V_n^p)}
- \sum_{p=1}^d \sum_{  {\scriptstyle I\subseteq \{1,\ldots,d\}:} \atop {\scriptstyle \|I\|=p}   } \prod_{q\in I}  \mathbb{I}_{ \bigcup\limits_{i=0}^j ( 0\vee s_{i,j-1} - C\epsilon,  s_{i,j} + C\epsilon   ] }^{(V_n^p)}
\end{align}
\end{subequations}
for all $k$ sufficiently large.
Using both inequalities
and that
$$
\mathbb{I}_{ ( 0, 1-Z_{j}^{N_k}(t_n^{N_k,\lambda -}) ] \bigcup\, (1 - Z_{j+1}^{N_k}(t_n^{N_k,\lambda -}) ,1]}^{(V_n^p)}
\times
\mathbb{I}_{ ( 0, 1-Z_{j}^{N_k}(t_n^{N_k,\lambda -}) ] }^{(V_n^p)}
=
\mathbb{I}_{ (1 - Z_{j+1}^{N_k}(t_n^{N_k,\lambda -}) ,1]}^{(V_n^p)}, 
$$
for the term \eqref{ascjh87tf} we obtain
\begin{align*}
&R_0(t) - R_1(t) - \lim_{\epsilon\to 0} \frac{1}{\epsilon} \times O(\epsilon) \epsilon \\
&\le  \lim_{\epsilon\to 0} \frac{1}{\epsilon}\lim_{k\to\infty}
 \frac{1}{N_k} \sum_{n=\mathcal{N}_\lambda(N_k t)+1}^{\mathcal{N}_\lambda(N_k (t+\epsilon))} 
\mathbf{1}_{ \left\{R_{0}^{N_k}(t_n^{N_k,\lambda-}) = 0 \right\} } 
\left (  F_{0,1,n}^{N_k} + F_{1,1,n}^{N_k}  \right) 
\prod\limits_{p=1}^d \mathbb{I}_{ ( 1-Z_1^{N_k}(t_n^{N_k,\lambda-}), 1 ] }^{(V_n^p)} \\
&\le 
R_0(t) - R_1(t).
\end{align*}
where the $O(\epsilon) \epsilon$ term is obtained by applying Lemma~\ref{GC_F} to the terms in the double sum of \eqref{sjc8asc9h}, which gives 
%
\begin{align*}
\lim_{k\to\infty}
 \frac{1}{N_k} \sum_{n=\mathcal{N}_\lambda(N_k t)+1}^{\mathcal{N}_\lambda(N_k (t+\epsilon))} 
\sum_{p=1}^d \sum_{  {\scriptstyle I\subseteq \{1,\ldots,d\}:} \atop {\scriptstyle \|I\|=p}   } \prod_{q\in I}  \mathbb{I}_{ \bigcup\limits_{i=0}^j ( 0\vee s_{i,j-1} - \epsilon C,  s_{i,j} + \epsilon C   ] }^{(V_n^p)}
= 
\sum_{p=1}^d \sum_{  {\scriptstyle I\subseteq \{1,\ldots,d\}:} \atop {\scriptstyle \|I\|=p}   }
(2\epsilon C)^p \times \lambda \epsilon = O(\epsilon) \epsilon.
\end{align*}
Thus, when $\overline{x}_{0,0}(t)+\overline{x}_{\cdot,1}(t)=0$ and $t$ is a point of differentiability, we obtain
\begin{equation*}
\lim_{\epsilon\downarrow 0} \frac{1}{\epsilon}\lim_{k\to\infty}
 \frac{1}{N_k} \sum_{n=\mathcal{N}_\lambda(N_k t)+1}^{\mathcal{N}_\lambda(N_k (t+\epsilon))} 
\mathbf{1}_{ \left\{R_{0}^{N_k}(t_n^{N_k,\lambda-}) = 0 \right\} } 
\left (  F_{0,1,n}^{N_k} + F_{1,1,n}^{N_k}  \right) 
\prod\limits_{p=1}^d \mathbb{I}_{ ( 1-Z_1^{N_k}(t_n^{N_k,\lambda-}), 1 ] }^{(V_n^p)}
=
R_0(t) - R_1(t)
\end{equation*}
and substituting this in  \eqref{ascjh87tfasc} we obtain
\begin{align*}
 0= \lim_{\epsilon\to 0} \frac{1}{\epsilon}\lim_{k\to\infty}  X_{\cdot,1}^{N_k}(t+\epsilon)-X_{\cdot,1}^{N_k}(t) =  \lambda d \overline{x}_{1,\cdot}(t)   + \lambda -  2R_0(t)   +  R_1(t),
\end{align*}
which implies that $R_1(t)$ exists and furthermore that (since $R_1(t)$ is necessarily non-negative by definition)
is given by
\begin{align*}
  R_1(t) &  = 0\vee \left( 2 R_0(t) -  \lambda d \overline{x}_{1,\cdot}(t)   -\lambda  \right)
\mathbf{1}_{\{\overline{x}_{0,0}(t)+\overline{x}_{\cdot,1}(t)=0\}}\\
& = 0\vee \lambda \left( 1 -2 d\overline{x}_{0,\cdot}(t)) -  d \overline{x}_{1,\cdot}(t)   \right)
\mathbf{1}_{\{\overline{x}_{0,0}(t)+\overline{x}_{\cdot,1}(t)=0\}}
\end{align*}
where the last equation follows by substituting  the expression of $R_0(t)$ given in Lemma~\ref{prop_R0}.
This proves \eqref{eq:sdjviuy7r1} when $j=1$.

If $j>1$, the same argument applies again. First, 
taking summations over the sample path expressions~\eqref{eq:shucyu} and~\eqref{eq:shucyuascs}, we obtain
\begin{subequations}
\nonumber
\begin{align}
X_{\cdot,j}^N(t)  = & X_{\cdot,j}^N(0 )  + \frac{1}{N} \sum_{n=1}^{\mathcal{N}_\lambda(N t)}
\sum_{p=1}^d  \mathbb{I}_{ ( S_{j,j}^N(t_n^{N,\lambda-})  , S_{j,I}^N(t_n^{N,\lambda-})  ] }^{(V_n^p)}\\
\\
&- \sum_{i=0}^{j-1}   \frac{1}{N} \sum_{n=1}^{\mathcal{N}_\lambda(N t)}
\sum_{p=1}^d  \mathbb{I}_{ ( S_{i,j}^N(t_n^{N,\lambda-})-X_{i,j}^N(t_n^{N,1-}),  S_{i,j}^N(t_n^{N,\lambda-})   ] }^{(V_n^p)}\\
%
%
& + \frac{1}{N} \sum_{n=1}^{\mathcal{N}_\lambda(N t)}  \mathbf{1}_{ \left\{R_{j-2}^N(t_n^{N,\lambda-}) = 0 \right\} } 
\left(   \sum_{i=0}^{j-1} F_{i,j-1,n}^N  \right) 
\prod\limits_{p=1}^d \mathbb{I}_{ ( 1-Z_{j-1}^N(t_n^{N,\lambda-}), 1 ] }^{(V_n^p)} \\
& - \frac{1}{N} \sum_{n=1}^{\mathcal{N}_\lambda(N t)}  \mathbf{1}_{ \left\{R_{j-1}^N(t_n^{N,\lambda-}) = 0 \right\} } 
\left( \sum_{i=0}^j F_{i,j,n}^N \right) 
\prod\limits_{p=1}^d \mathbb{I}_{ ( 1-Z_{j}^N(t_n^{N,\lambda-}), 1 ] }^{(V_n^p)}.
\end{align}
\end{subequations}
Then, using \eqref{ascuja7sc8as}
and noting that
the term in \eqref{fjw8tg8irf} can be bounded as in \eqref{jvr7h6tg}
we obtain
\begin{align*}
0=& \dot{\overline{x}}_{\cdot,j}(t)= \lim_{\epsilon\to 0}\frac{1}{\epsilon}\lim_{k\to\infty}  X_{\cdot,j}^{N_k}(t+\epsilon)-X_{\cdot,j}^{N_k}(t) \\
=&  \lambda d \overline{x}_{j,\cdot}(t)   + (R_{j-2}(t) - R_{j-1}(t) ) - (R_{j-1}(t)  - R_{j}(t) ),
\end{align*}
provided that $\sum_{j'=0}^j \overline{x}_{\cdot,j'}(t)=0$ and $t$ is a point of differentiability.
This condition inductively proves the existence of $R_{j}(t) )$ and implies that
\begin{align}
\nonumber
 R_{j}(t) = & 0\vee 
\left(
2R_{j-1}(t) - R_{j-2}(t) - \lambda d \overline{x}_{j,\cdot} 
\right)
\mathbf{1}_{\{\sum_{j'=0}^j \overline{x}_{\cdot,j'}(t)=0\}}.
\end{align}
Since we have already obtained an expression for $R_0$ and $R_1$, we can derive an expression for $R_{2}$ and so forth iteratively for all $j$.
By iteratively substituting in previous equation the expressions of $R_{j-1}$ and $R_{j-2}$, we obtain \eqref{eq:sdjviuy7r1}.

\section{Proof of Lemma~\ref{lemma:iv8fd}}
As shown at the beginning of the proof of Lemma~\ref{prop_R0}, if $\sum_{i=0}^{j} \overline{x}_{\cdot,i}(t)>0$, then 
$R_{j}^{N_k}(t_n^{N_k,\lambda-})>0$ for all~$k$ sufficiently large, and therefore $\Gamma_{i,j}=0$.
Thus, let us assume in the following that $\sum_{i=0}^{j} \overline{x}_{\cdot,i}(t)=0$.
We give a proof when $i\le j$; when $i =j+1$, the proof uses the same argument and is omitted.

So far we have assumed that $\omega \in \mathcal{C}$ was fixed but for now let us explicit the dependence on $\omega$ and treat quantities $\overline{x}(t)$ and $X^N(t)$ as random variables.

Let
\begin{align}
\label{ascu876tf}
\Gamma_{i,j}^{\epsilon,N}(t)
\bydef
\frac{1}{\epsilon N}
\sum_{n=\mathcal{N}_\lambda(N t)+1}^{\mathcal{N}_\lambda(N (t+\epsilon))} 
\underbrace{
\mathbf{1}_{ \left\{R_{j}^{N}(t_n^{N,\lambda-}) = 0 \right\} }
F_{i,j+1,n}^N
\prod_{p=1}^d \mathbb{I}_{(1 - Z_{j+1}^{N}(t_n^{N,\lambda -}) ,1]}^{(V_n^p)}
}_{\bydef Y_{n,N}}
\end{align}
For all $n$, the random variable
$Y_{n,N}$
is  $\mathcal{F}_n$-measurable where
$\mathcal{F}_n \bydef \{X^{N}(t_n^{N,\lambda-}), V_n^1,\ldots,V_n^d, U_n,  W_n\}$, and 
\begin{align}
\label{asjc7h6ftc}
\mathbb{E}[Y_{n,N}| \mathcal{F}_n\setminus W_n] 
=
\mathbf{1}_{ \left\{R_{j}^{N}(t_n^{N,\lambda-}) = 0 \right\} }
%
%
\frac
{ X_{i,j+1}^N(t_n^{N,\lambda-}) - M_{i,j+1,n}^N }
{X_{\cdot,j+1}^N(t_n^{N,\lambda-}) + \overline{M}_{j+1,n}^N - \underline{M}_{j+1,n}^N}
\prod_{p=1}^d \mathbb{I}_{(1 - Z_{j+1}^{N}(t_n^{N,\lambda -}) ,1]}^{(V_n^p)}
\end{align}
because we recall that the 0-1 random variable $F_{i,j,n}^N$ is one if and only if $\big( X_{\cdot,j}^N(t_n^{N,\lambda-}) + \overline{M}_{j,n}^N - \underline{M}_{j,n}^N\big) \, W_n \in  \big( \sum_{k=0}^{i-1} X_{k,j}^N(t_n^{N,\lambda-}) -M_{k,j,n}^N,  \sum_{k=0}^{i} X_{k,j}^N(t_n^{N,\lambda-}) - M_{k,j,n}^N  \big] $, by definition \eqref{def:F_ijn}, with~$W_n$ uniform over [0,1]; the set $\mathcal{F}_n\setminus W_n$ denotes the set $\mathcal{F}_n$ with $W_n$ removed.

Let $Z_{n,N}\bydef Y_{n,N} - \mathbb{E}[Y_{n,N}| \mathcal{F}_n\setminus W_n]$.
Then, $\mathbb{E}[Z_{n,N}|\mathcal{F}_n\setminus W_n]=0$ and $|Z_{n,N}|\le 2$, and applying the Azuma--Hoeffding inequality, we get
\begin{equation}
\mathbb{P} \left( \frac{1}{N} \left|\sum_{n=1}^N Z_{n,N}  \right| >\delta \right) \le 2 \exp\left( -\frac{(N\delta)^2}{ 8 N} \right)
\end{equation}
for any $\delta>0$.
Since $\sum_N \exp\left( -{N \delta^2/ 8} \right)<\infty$, an application of the Borel--Cantelli lemma shows that
$\frac{1}{N} \sum_{n=1}^N Z_{n,N} \to 0$ almost surely.
In particular,
\begin{equation}
\label{asjc8as9}
\lim_{N \to\infty} \Gamma_{i,j}^{\epsilon,N}(t) 
-
\frac{1}{N\epsilon}
\sum_{n=\mathcal{N}_\lambda(N t)+1}^{\mathcal{N}_\lambda(N (t+\epsilon))} 
\mathbb{E}[Y_{n,N}| \mathcal{F}_n\setminus W_n]  = 0
\end{equation}
almost surely.

Now, we fix $\omega\in\mathcal{C}$ and use \eqref{asjc7h6ftc} and Lemma~\ref{lemma:X_close_x} to obtain that for any $\epsilon>0$
\begin{align*}
& \lim_{k\to\infty}\frac{1}{\epsilon N_k}
\sum_{n=\mathcal{N}_\lambda(N_k t)+1}^{\mathcal{N}_\lambda(N_k (t+\epsilon))} 
 \mathbb{E}[Y_{n,N_k}| \mathcal{F}_n\setminus W_n] \\
 & \le
 \lim_{N_k\to\infty}\frac{1}{\epsilon N_k}
\sum_{n=\mathcal{N}_\lambda(N_k t)+1}^{\mathcal{N}_\lambda(N_k (t+\epsilon))} 
\mathbf{1}_{ \left\{R_{j}^{N_k}(t_n^{N_k,\lambda-}) = 0 \right\} }
%
%
\frac
{ \overline{x}_{i,j+1}(t) + C\epsilon }
{ \overline{x}_{\cdot,j+1}(t) - C\epsilon}
\,\prod_{p=1}^d \mathbb{I}_{(1 - Z_{j+1}^{N_k}(t_n^{N_k,\lambda -}) ,1]}^{(V_n^p)}
\end{align*}
Replacing $\epsilon$ by $-\epsilon$ in the last fraction term,
the previous inequality can be reversed and letting $\epsilon\downarrow 0$,
we obtain
\begin{align}
\label{vjuygedffheryf6723ws}
\lim_{\epsilon\downarrow 0} \lim_{k\to\infty}\frac{1}{\epsilon N_k}
\sum_{n=\mathcal{N}_\lambda(N_k t)+1}^{\mathcal{N}_\lambda(N_k (t+\epsilon))} 
 \mathbb{E}[Y_{n,N_k}| \mathcal{F}_n\setminus W_n]
 =
R_j(t)%
\frac
{ \overline{x}_{i,j+1}(t)}
{ \overline{x}_{\cdot,j+1}(t) }.
\end{align}
Finally, \eqref{asjc8as9} and \eqref{vjuygedffheryf6723ws} give \eqref{asjc78dcckiujff}.

\section*{Acknowledgments}
The authors would like to thank Sem Borst, Bruno Gaujal and the referees for their valuable comments and remarks.


\bibliographystyle{abbrv} 

\end{document}